\documentclass[a4paper,12pt]{article}
\usepackage[utf8]{inputenc}
\usepackage[left=2cm, right=2cm, top=2cm, bottom=2cm]{geometry}
\usepackage{amssymb}

\usepackage{bm}
\usepackage[normalem]{ulem}
\usepackage[justification=centering]{caption}
\usepackage{subfig}
\usepackage{optidef}
\usepackage{epstopdf}

\usepackage{stanli}
\usepackage{adjustbox}
\usepackage{graphicx}
\usepackage{multirow}
\usepackage{multicol}
\usepackage{hhline}
\usepackage{adjustbox}
\usepackage{amssymb}
\usepackage{amsthm}
\usepackage{thmtools}
\usepackage{mathdots}
\usepackage{enumitem}

\newtheorem{theorem}{Theorem}
\newtheorem{proposition}[theorem]{Proposition}
\newtheorem{conjecture}[theorem]{Conjecture}
\newtheorem{lemma}[theorem]{Lemma}
\newtheorem{remark}[theorem]{Remark}

\newtheorem{definition}[theorem]{Definition}

\usepackage{xcolor}
\usepackage{hyperref}
  \hypersetup{
  colorlinks,
  citecolor=blue,
  linkcolor=blue,
  urlcolor=blue}

\newcommand{\R}{\mathbb{R}}
\newcommand{\C}{\mathbb{C}}
\newcommand{\N}{\mathbb{N}}
\newcommand{\B}{\mathcal{B}}
\newcommand{\Pc}{\mathcal{P}}
\newcommand{\cJ}{\mathcal{J}}

\newcommand{\Gb}{\mathcal{G}}
\newcommand{\xb}{\mathbf{x}}
\newcommand{\wb}{\mathbf{w}}
\newcommand{\zb}{\mathbf{z}}
\DeclareMathOperator{\rank}{rank}
\DeclareMathOperator{\LT}{LT}
\DeclareMathOperator{\kernel}{Ker}
\DeclareMathOperator{\mom}{mom}
\DeclareMathOperator{\sos}{sos}
\DeclareMathOperator{\lex}{lex}

\def\hone{x}
\def\htwo{w}

\def\rfig#1{Fig.~\ref{#1}}

\def\rth#1{Theorem~\ref{#1}}

\def\rpr#1{Proposition~\ref{#1}}

\def\rle#1{Lemma~\ref{#1}}
\def\rdef#1{Definition~\ref{#1}}

\def\rsub#1{Subsection~\ref{#1}}

\def\He{{\rm He}}

\def\diag{{\rm diag}}
\def\trace{{\rm trace}}
\def\rank{{\rm rank}}

\def\Ra{\Rightarrow}

\def\PSD{{\mathcal{PSD}}}
\def\NN{{\mathcal{NN}}}
\def\COP{{\mathcal{COP}}}

\def\bPi{\boldsymbol{\Pi}}

\def\bfPi{{\mathbf \Pi}}
\def\bfPis{{\bfPi^\star}}

\def\clIu{\clI_\mathrm{u}}
\def\clIl{\clI_\mathrm{l}}

\def\rfig#1{Fig.~\ref{#1}}

\def\rth#1{Theorem~\ref{#1}}

\def\rpr#1{Proposition~\ref{#1}}

\def\rle#1{Lemma~\ref{#1}}
\def\rdef#1{Definition~\ref{#1}}

\def\rsub#1{Subsection~\ref{#1}}

\def\He{{\rm He}}

\def\diag{{\rm diag}}
\def\trace{{\rm trace}}
\def\rank{{\rm rank}}

\def\Ra{\Rightarrow}

\def\PSD{{\mathcal{PSD}}}
\def\NN{{\mathcal{NN}}}
\def\COP{{\mathcal{COP}}}

\def\bPi{\boldsymbol{\Pi}}

\def\bfPi{{\mathbf \Pi}}
\def\bfPis{{\bfPi^\star}}

\def\clIu{\clI_\mathrm{u}}
\def\clIl{\clI_\mathrm{l}}

\def\bbD{{\mathbb D}}

\def\bbH{{\mathbb H}}

\def\bbR{{\mathbb R}}
\def\bbS{{\mathbb S}}

\def\clH{{\mathcal H}}
\def\clI{{\mathcal I}}
\def\clJ{{\mathcal J}}

\title{LMI hierarchies for stability analysis of ReLU feedback systems}
\begin{document}

\author{
Victor Magron$^1$ \and 
Yoshio Ebihara$^2$ \and 
Shingo Nishinaka$^2$ \and 
Dimitri Peaucelle$^1$ \and
Rin Saeki$^2$ \and  
Tsuyoshi Yuno$^2$ \and 
Sophie Tarbouriech$^1$ 
}
\footnotetext[1]{Universit\'{e} de Toulouse; LAAS-CNRS, 7 avenue du colonel Roche, F-31400 Toulouse, France}
\footnotetext[2]{Graduate School of Information Science \& Electrical Engineering, Kyushu University, 744 Motooka, Nishi-ku, Fukuoka 819-0395, Japan}
\date{\today}

\maketitle

\begin{abstract}
We consider the stability analysis of feedback systems with rectified linear unit (ReLU) activations, and model this problem with polynomial optimization. 
Stability can be certified by means of copositive multipliers in the framework of integral quadratic constraints. 
Based on a duality argument, we show how to certify instability by considering a complete hierarchy of linear matrix inequalities. 
This hierarchy is obtained by leveraging the specific equality constraints arising from the ReLU encoding. 
We illustrate the effectiveness of the proposed approach through several numerical examples. 
\end{abstract}


\section{Introduction}
\label{sec:intro}

Control theoretic approaches for the analysis and synthesis of neural networks (NNs)  
(see, e.g., \cite{Fazlyab_IEEE2022,Raghunathan_NIPS2018,Revay_LCSS2021,Gronqvist_MTNS2022,Ebihara_EJC2021,Ebihara_arXiv2024a})
and of dynamical systems driven by NNs \cite{Yin_IEEE2022,Scherer_IEEEMag2022} 
have recently attracted great attention.  
In this study, we are particularly interested in dynamical NNs 
such as recurrent NNs (RNNs) and 
NN-driven control systems where the activation functions in NNs are
rectified linear units (ReLUs).  
Such systems can be regarded as a special class of 
Lurie systems for which effective 
analysis and synthesis conditions have been developed, see, e.g.,  
\cite{Tarbouriech_2011}.  
The generality of the approaches for Lurie systems comes 
at the expense of strong conservativeness,   
and studies dedicated to 
ReLU nonlinearity are relatively scarce.  
In this paper, we aim to derive analysis conditions
that are specific to ReLU nonlinearity.   

We focus on the semialgebraic set representation
that accurately characterizes the input-output properties of ReLUs. 
Such representation has been already successfully combined with semidefinite programming in 
\cite{Raghunathan_NIPS2018,
chen2020semialgebraic} to analyze the robustness of ReLU  networks, and in \cite{Groff_CDC2019,korda2022stability} to certify the stability of dynamical systems controlled by ReLU networks. 

In \cite{Ebihara_arXiv2024a}, the authors rely on the same representation and introduce the so called \textit{copositive multiplier} 
that captures the input-output properties in quadratic form.   
By using this novel multiplier 
in the general framework of integral quadratic constraint (IQC) theory 
(\cite{Megretski_IEEE1997,Scherer_IEEEMag2022}), 
one can derive a (primal) linear matrix inequality (LMI) condition for the stability analysis of feedback systems with ReLU nonlinearities.  
However, this primal LMI is conservative in general.  
Therefore, if the primal LMI turns out to be infeasible by numerical computation, one cannot conclude anything about the system's stability.   
This motivates us to focus on the dual LMI whose feasibility is guaranteed by the theorems of alternatives.  
By closely investigating the structure of the dual solution, we obtain a rank condition on the dual variable under which we can conclude that the feedback system is not globally asymptotically stable.  
More precisely, from the dual solution satisfying the rank condition, one can extract an initial condition from which  the state trajectory does not converge to the origin.  
In robust stability analysis of linear dynamical systems affected by parametric uncertainties, it is well known that dual LMIs are quite useful in extracting worst-case parameters that destabilize the systems, see, e.g., 
\cite{Scherer_SIAM2005,Masubuchi_SCL2009,Ebihara_IEEE2009}.  
In the present paper, the dual LMI works effectively 
in extracting the worst-case initial condition 
for nonlinear systems in the above sense.  
We also show that the dual LMI can be interpreted 
as a relaxation of a system of polynomial (in)equalities 
that characterizes the existence of the worst-case initial conditions.  
This enables the construction of a hierarchy of higher order LMIs that progressively improves the ability to detect instability.

\paragraph{Outline and novelty with respect to \cite{relufeedback}}
This paper is an extended version of the article \cite{relufeedback} published in the proceedings of the MICNON conference, in which we derived LMIs for (in)stability certification of ReLU feedback systems. \\
The contents of Section  \ref{sec:prelim}, Section \ref{sec:stability} and Section \ref{sec:higher} already appear in \cite{relufeedback}. 
Section \ref{sec:prelim} sets up the problem formulation and introduces the semialgebraic set representation of ReLU nonlinearities, which forms the foundation for our analysis. 
Section \ref{sec:stability} presents the primal and dual LMI conditions for stability and instability certification, respectively, including a rank-one condition under which the dual LMI guarantees instability. 
Section \ref{sec:higher} interprets the dual LMI as the first-order relaxation of a nonconvex polynomial feasibility problem and introduces a hierarchy of higher-order LMI relaxations based on block-Hankel structures, designed to improve the likelihood of detecting instability.
We illustrate the effectiveness of the first and second order relaxations by several numerical examples in Section \ref{sec:num} and Section \ref{sub:num}, respectively. \\
In this extended version, we provide a more comprehensive theoretical foundation for the instability analysis. 
Specifically, we introduce in Section \ref{sec:setting} a new moment-sums of squares hierarchy of LMIs based on the specific structure exploitation of the equality constraints, in the spirit of \cite{gradsos}. 
This hierarchy complements the dual LMI hierarchy presented in Section \ref{sec:higher}, and offers an alternative pathway to instability certification with formal convergence guarantees. 

\paragraph{Notation} 
The set of $n\times m$ real matrices is denoted by $\bbR^{n\times m}$, and
the set of $n\times m$ entrywise nonnegative matrices is denoted
by $\bbR_+^{n\times m}$. 
For a matrix $A$, we also write $A\geq 0$ to denote that  
$A$ is entrywise nonnegative. 
We denote the set of $n\times n$ real symmetric, positive semidefinite (definite), and Hurwitz stable matrices 
by $\bbS^n$,  $\bbS_{+}^n$ ($\bbS_{++}^n$), and $\bbH^n$, respectively.  
For $A\in\bbS^n$, we also write $A\succ 0\ (A\prec 0)$ to
denote that $A$ is positive (negative) definite.  
For $A\in\bbR^{n\times n}$, we define $\He\{A\}:=A+A^T$ and 
$\diag(A)\in\bbR^n$ where $(\diag(A))_i:=A_{i,i}\ (i=1,\cdots,n)$.  
For $A\in\bbR^{n\times n}$ and $B\in\bbR^{n\times m}$,
$(\ast)^TAB$ is a shorthand notation of $B^TAB$.  
We denote by $\bbD^n\subset\bbR^{n\times n}$ the set of diagonal matrices.  
For $v\in\bbR^{n}$, $\|v\|$ stands for the standard Euclidean norm.    
For a static (possibly nonlinear) operator $\Psi:\ \bbR^n\mapsto \bbR^m$, 
we denote by $\|\Psi\|$ its $l_2$-induced norm.  
We define the positive semidefinite cone $\PSD^{n}\subset \bbS^n$, 
the copositive cone $\COP^{n}\subset \bbS^n$, and the nonnegative cone
$\NN^{n}\subset \bbS^n$ as follows:
\[
 \begin{array}{@{}l}
  \PSD^n:=\{P\in\bbS^n:\ x^TPx\geq 0\ \forall x\in\bbR^n\},\\
  \COP^n:=\{P\in\bbS^n:\ x^TPx\geq 0\ \forall x\in\bbR_{+}^n\},\\
  \NN^n:= \{P\in\bbS^n:\ P\geq 0\}.  
 \end{array}
\]
We can readily see that $\PSD^n\subset \PSD^n+\NN^n\subset \COP^n$,
where ``$+$'' here stands for the Minkowski sum.

\section{Problem setting and preliminary results}
\label{sec:prelim}
\subsection{Problem statement}
\begin{figure}[b]
\begin{center}
\if{
\begin{picture}(3.5,2.5)(0,0)
\put(0,2){\vector(1,0){1}}
\put(1,1.5){\framebox(1.5,1){$\Phi$}}
\put(2.5,2){\line(1,0){1}}
\put(3.5,2){\line(0,-1){1.5}}
\put(3.5,0.5){\vector(-1,0){1}}
\put(3,0.3){\makebox(0,0)[t]{$w$}}
\put(1,0){\framebox(1.5,1){$G$}}
\put(1,0.5){\line(-1,0){1}}
\put(0.5,0.3){\makebox(0,0)[t]{$z$}}
\put(0,0.5){\line(0,1){1.5}}
\end{picture}
}\fi
\begin{tikzpicture}[auto, node distance=1.5cm and 2cm, >=latex']
  \node [draw, minimum width=1.5cm, minimum height=1cm] (phi) {$\Phi$};
  \node [draw, below=of phi, minimum width=1.5cm, minimum height=1cm] (G) {$\textbf{G}$};

  \draw[<-] (phi.west) -- ++(-1,0) -- ++(0,-2.5) --  (G.west);
  \draw[<-] (G.east) -- ++(1,0) -- ++(0,2.5) -- (phi.east);

\node at (1.2,-2.8) {$\wb$};
\node at (-1.2,-2.8) {$\zb$};

\end{tikzpicture}
\caption{Nonlinear Feedback System $\Sigma$.}
\label{fig:Sigma}
\end{center}
\end{figure}
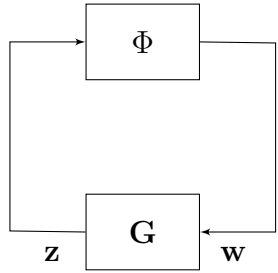

Let us consider the feedback system $\Sigma$ shown in \rfig{fig:Sigma}.  
Here, $\textbf{G}$ is a linear system described by 
\begin{equation}
 \mathbf{G}:\ 
 \left\{
 \begin{array}{lcl}
  \dot \xb(t) & = & A \xb(t)+B \wb(t), \\
  \zb(t) & = & C \xb(t)+D\wb(t)  
 \end{array}\right.  
\label{eq:G}
\end{equation}
where 
$A\in\bbR^{n\times n}\cap \bbH^n$, 
$B\in\bbR^{n\times m}$, 
$C\in\bbR^{m\times n}$, and   
$D\in\bbR^{m\times m}$.   
Moreover, we have
\begin{equation}
 \wb(t) =\Phi(\zb(t)), 
\label{eq:Phi}
\end{equation}
where $\Phi:\ \bbR^m\mapsto \bbR^m$ 
stands for the static rectified linear unit (ReLU) nonlinearity given by
\begin{equation}
\begin{array}{@{}l}
 \Phi(q)=\left[\ \phi(q_1)\ \cdots\ \phi(q_m)\ \right]^T,\\
 \phi: \bbR\to \bbR,\quad
 \phi(\eta)=
  \left\{
   \begin{array}{cc}
    \eta & (\eta\ge 0),  \\
    0 & (\eta< 0).  \\
   \end{array}
  \right.    
\end{array}  
\label{eq:ReLU}
\end{equation}
\begin{remark}
We note that ReLU $\Phi:\bbR^m\mapsto\bbR^m$
satisfies
\begin{equation}
\Phi(\alpha q) = \alpha \Phi(q)\ \forall (\alpha,q)\in\bbR_+\times \bbR^m.  
\label{eq:co-linear}
\end{equation}
In addition, for $p,q\in\bbR^m$, 
it is shown in \cite{Raghunathan_NIPS2018} and \cite{Groff_CDC2019} that 
$p=\Phi(q)$ holds if and only if 
\begin{equation}
p-q\ge 0,\ p\ge 0,\ (p-q)\odot p =0  
\label{eq:ReLU_alg}
\end{equation}
where $\odot$ stands for the Hadamard product (element-wise product).  
Namely, the input-output properties of the ReLU $\Phi$ can be
accurately characterized by the semialgebraic set \eqref{eq:ReLU_alg}.  
\end{remark}
It can be easily seen that $\|\Phi\|=1$.  
With this fact in mind, for the well-posedness of the feedback system $\Sigma$, 
we assume $\|D\|<1$. 
Under this assumption, we see that the map 
$(I-D\Phi):\ \bbR^m\mapsto \bbR^m$  is invertible and hence 
the inverse map $(I-D\Phi)^{-1}:\ \bbR^m\mapsto \bbR^m$ is well defined.  
Moreover, the next results follow.  
\begin{lemma}
Suppose $\|D\|<1$.   
Then, we can represent the dynamics of the system $\Sigma$ by
\begin{equation}
\Sigma:\ \dot \xb(t)=A\xb(t)+B\Phi\circ(I-D\Phi)^{-1}(C\xb(t)).  
\label{eq:Sigma} 
\end{equation}
In addition, we have
\begin{equation}
 (I-D\Phi)^{-1}(\alpha q) = \alpha (I-D\Phi)^{-1}(q)\ \forall (\alpha,q)\in\bbR_+\times \bbR^m.  
\label{eq:co-linear-invmap}
\end{equation}
\label{le:Sigma}
\end{lemma}
\begin{proofof}{\rle{le:Sigma}}
It is obvious that \eqref{eq:Sigma} since $\zb(t) = C \xb(t) + D \Phi(\zb(t)) $  implies $(I - D \Phi)^{-1} (C \xb(t)) = (I - D \Phi)^{-1})(I - D \Phi) \zb(t) = \zb(t)$. 
For the proof of \eqref{eq:co-linear-invmap}, 
for given $q\in \bbR^m$,   
we first define
\[
 p:=(I-D\Phi)^{-1}(q).  
\]
This implies $q=(I-D\Phi)(p)$.  Then, for given $\alpha\in\bbR_+$, we have
from \eqref{eq:co-linear} that 
\[
 \alpha q = \alpha (I-D\Phi)(p)=  (I-D\Phi)(\alpha p).  
\]
This clearly shows that \eqref{eq:co-linear-invmap} holds.  
\end{proofof}
\subsection{Stability analysis}
This paper addresses the stability analysis problem of the system $\Sigma$.  
The definition of ``stability'' is now recalled.  

\begin{definition}[\cite{Khalil_2002}]
The equilibrium point $x=0$ (i.e., the origin) of the feedback system $\Sigma$ 
given by \eqref{eq:Sigma} is said to be locally asymptotically stable if the
conditions (i) and (ii) given below are satisfied, 
and is said to be globally asymptotically stable if the
conditions (i) and (iii) given below are satisfied, where
\begin{itemize}
 \item[(i)] for each $\varepsilon>0$, there exists $\delta=\delta(\varepsilon)>0$ such that
	    \[
	     \|\xb(0)\|<\delta\ \Rightarrow\ \|\xb(t)\|<\varepsilon\ (\forall t\ge 0);
	    \]
 \item[(ii)] there exist $\delta>0$ such that 
	    \[
	     \|\xb(0)\|<\delta\ \Rightarrow\ \lim_{t\to \infty}\xb(t)=0;
	    \]
 \item[(iii)] $\lim_{t\to \infty}\xb(t)=0$ holds for any $\xb(0)$.  
\end{itemize}
\label{def:stab}
\end{definition}

For general nonlinear feedback systems, it is of course true that 
there is a clear distinction between the local and global asymptotic stability.  
However, for the nonlinear feedback system $\Sigma$ with ReLU nonlinearities, 
they are equivalent as shown in the next theorem.  
\begin{theorem}
The nonlinear feedback system $\Sigma$ given by 
\eqref{eq:Sigma} is globally asymptotically stable if and only if 
it is locally asymptotically stable.
\label{th:gstab}
\end{theorem}
\begin{proofof}{\rth{th:gstab}}
Let us define the state trajectory of \eqref{eq:Sigma} 
for the initial state $\xb(0)=\xi$ by $\xb_\xi(t)\ (t\ge 0)$.  
Then, from the definition of the 
local and global asymptotic stability given in 
\rdef{def:stab}, 
we see that \rth{th:gstab} is verified 
by showing that $\xb_{\alpha \xi}(t)=\alpha \xb_{\xi}(t)$ holds for any 
$\alpha \in\bbR_+$.  
This can be readily validated since
\[
\begin{aligned}
 \frac{d}{dt} (\alpha \xb_{\xi}(t))&=\alpha \dot \xb_{\xi}(t)\\
&=\alpha \left(A\xb_\xi(t)+B\Phi\circ(I-D\Phi)^{-1}(C \xb_\xi(t))\right)\\
&=A(\alpha \xb_\xi(t))+B\Phi\circ(I-D\Phi)^{-1}(C(\alpha \xb_\xi(t)))  
\end{aligned}
\]
where we used \eqref{eq:co-linear} and \eqref{eq:co-linear-invmap}.  
This clearly shows that $\xb_{\alpha \xi}(t)=\alpha \xb_{\xi}(t)$ holds and hence
the proof is completed.  
\end{proofof}
According to the proof of Theorem \ref{th:gstab}, if a trajectory starting at $\xb(0)=\xi$ does not converge to the equilibrium point, this is also the case for all trajectories from the open ray $\xb(0)=\alpha \xi$ with $\alpha > 0$. Any such $\xi$ proves instability. 

On the basis of \rth{th:gstab}, we investigate the global asymptotic stability analysis
of the feedback system $\Sigma$.  
For simplicity, we say that the feedback system $\Sigma$ is stable if it is 
globally asymptotically stable.  
The next proposition forms the basis for the stability analysis 
using integral quadratic constraint (IQC) theory.  
\begin{proposition}[\cite{Megretski_IEEE1997}]\ \\
Let us define $\bfPis\subset\bbS^{2m}$ by 
\begin{equation}
\begin{array}{@{}l}
\bfPis:=
\left\{
\Pi\in\bbS^{2m}:\ 
\left[
 \begin{array}{c}
  \zeta\\
  \xi\\
 \end{array}
\right]^T\Pi
\left[
 \begin{array}{c}
  \zeta\\
  \xi\\
 \end{array}
\right]\ge 0 \right. \\
\left. \hspace*{10mm} \forall
\left[
 \begin{array}{c}
  \zeta\\
  \xi\\
 \end{array}
\right]\in\bbR^{2m}\ \mathrm{s.t.}\ 
\xi = \Phi(\zeta)
\right\}.   
\end{array}
\label{eq:multi}
\end{equation}
The system $\Sigma$ is stable if there exist 
$P \succ 0$ and $\Pi\in\bfPis$ such that 
\begin{equation}
\begin{bmatrix} PA + A^TP & PB \\ B^TP & 0 \end{bmatrix} +
\begin{bmatrix} C & D \\ 0 & I_m \end{bmatrix}^T \Pi
\begin{bmatrix} C & D \\ 0 & I_m \end{bmatrix} \prec 0.   
\label{eq:stab}
\end{equation}
\label{pr:basic}
\end{proposition}
%

\subsection{Multipliers capturing I/O ReLU properties}

In the IQC-based stability condition \eqref{eq:stab}, 
it is of prime importance to employ a set of multipliers
$\bPi\subset\bPi^\star$
that is numerically tractable and 
captures the input-output properties of ReLU $\Phi$
as accurately as possible.  
On this issue, the next results have been obtained recently.  
\begin{proposition}[\cite{Ebihara_arXiv2024a}]\ \\
Let us define $\bPi_\COP,\bPi_\NN\subset\bbS^{2m}$ by
\begin{equation}
\scalebox{1.0}{$
\begin{array}{@{}l}
\bPi_\COP:=
\left\{\Pi\in\bbS:\ 
\Pi= E^T  
\left(
Q+\clJ(J)
\right)E,\right.\\ 
\left.\hspace*{30mm}
J\in\bbD^m,\ Q\in\COP^{2m}
\right\},\\
\bPi_\NN:=
\left\{\Pi\in\bbS:\ 
\Pi= E^T  
\left(
Q+\clJ(J)
\right)E,\right.\\ 
\left.\hspace*{30mm}
 J\in\bbD^m,\ Q\in\NN^{2m}
\right\},\\
E:=\left[
\begin{array}{ccc}
  -I_m & I_m \\
  0_{m,m} & I_m \\
\end{array}
\right],\ 
\clJ(J):=
\begin{bmatrix}
0_{m,m} & J \\
\ast & 0_{m,m} \\
\end{bmatrix}.  
\end{array}$}
\label{eq:PiCOP}
\end{equation}
Then we have $\bPi_\NN\subset\bPi_\COP\subset \bPi^\star$.  
\label{pr:PiCOP}
\end{proposition}

Propositions \ref{pr:basic} and \ref{pr:PiCOP} enable us
to derive a concrete LMI for the stability analysis of the system $\Sigma$;
see ``Primal LMI'' in \rth{th:pdLMI} in the next section.  

The validity of employing $Q\in\COP^{2m}$ in $\bfPi_\COP$
(or $Q\in\NN^{2m}$ in $\bfPi_\NN$) can be seen 
from the first and the second inequality constraints 
in the semialgebraic set representation
\eqref{eq:ReLU_alg}, and the validity of 
employing $J\in\bbD^{n}$ in $\bfPi_\COP$ and $\bfPi_\NN$
can be seen from the third equality constraint in  \eqref{eq:ReLU_alg}.  

As mentioned in \cite{Dur_2010}, the problem of determining whether
a given matrix is copositive or not  is a co-NP complete problem in general, and 
hence the set of multipliers $\bPi_\COP$ is numerically intractable.  
We therefore employ $\bPi_\NN$ in this paper.  
Since the ReLU $\Phi$ is a (repeated) slope-restricted nonlinearity, 
we can also employ known multipliers that are valid for this class of nonlinearities 
such as (static) O'Shea-Zames-Falb multipliers 
(\cite{O'Shea_IEEE1967,Zames_SIAM1968,Carrasco_EJC2016,Fetzler_IFAC2017})
and the multipliers proposed by  Fazlyab et al. (\cite{Fazlyab_IEEE2022}).    
It has been shown that $\bPi_\NN$ encompasses these known multipliers;
see \cite{Ebihara_arXiv2024a} for details.

\section{Primal-dual LMIs for (in)stability analysis}
\label{sec:stability}
\subsection{Theoretical results}
On the basis of Propositions~\ref{pr:basic} and \ref{pr:PiCOP}, 
we can obtain the next results.  
\begin{theorem}
Let us consider the primal and dual LMIs given respectively as follows:\\
\noindent
{\bf Primal LMI:}  
\begin{equation}
\mathrm{Find}\ P\in\bbS_{++}^n\ \mathrm{and}\ \Pi\in\bPi_\NN\ \mathrm{such\ that}\ \eqref{eq:stab}.  
\label{eq:pLMI}
\end{equation}
\noindent
{\bf Dual LMI:}  
\begin{equation}
\begin{array}{@{}l}
 \mathrm{Find}\ H=
\begin{bmatrix} H_{11} & H_{12} \\ H_{12}^T & H_{22} \end{bmatrix}
\in\bbS_+^{n+m}\setminus \{0\}\ \mathrm{such\ that}\\
\He\{AH_{11}+BH_{12}^T\}\succeq 0,\\ 
\begin{bmatrix} -C & -D+I_{m} \\ 0 & I_{m} \end{bmatrix} H
\begin{bmatrix} -C & -D+I_{m} \\ 0 & I_{m} \end{bmatrix}^T \ge 0,\\ 
\diag(-CH_{12}+(-D+I_m)H_{22}) =0.  
\end{array}
\label{eq:dLMI}
\end{equation}
The system $\Sigma$ is stable if the primal LMI \eqref{eq:pLMI} is feasible.  
In addition, the primal LMI \eqref{eq:pLMI} is infeasible if and only if the 
dual LMI \eqref{eq:dLMI} is feasible.  
\label{th:pdLMI}
\end{theorem}

The dual LMI \eqref{eq:dLMI} can be derived from the primal LMI \eqref{eq:pLMI}
by following the Lagrange duality theory 
for semidefinite programming problems (\cite{Scherer_EJC2006}).  
The last assertion in \rth{th:pdLMI} also follows from
theorems of alternative for LMIs; see, e.g., \cite{Scherer_EJC2006}.  

If the primal LMI \eqref{eq:pLMI} is feasible,
we can conclude that the system $\Sigma$ is stable.  
However, if the primal LMI \eqref{eq:pLMI} is infeasible,
we cannot conclude anything.  
This motivates us to focus on the dual LMI \eqref{eq:dLMI} that is
always feasible in this case as shown in \rth{th:pdLMI}.  
In fact, by paying attention to the structure of the dual solution, 
we can obtain the next theorem that is the first main result 
of this paper.  
\begin{theorem}
Suppose the dual LMI \eqref{eq:dLMI} is feasible 
with a dual solution $H\in\bbS_{+}^{n+m}\setminus\{0\}$ 
of $\rank(H)=1$.  Then, the system $\Sigma$ is not stable.  
Moreover, the following assertions hold:
\begin{itemize}
 \item[(i)] The full-rank factorization of $H$ is given of the form
	    \begin{equation}
	       H=\begin{bmatrix}\hone\\\htwo\end{bmatrix}
	       \begin{bmatrix}\hone\\\htwo\end{bmatrix}^T,\ 
	       \hone\in\bbR^n,\ \htwo\in \bbR_+^m. 
		\label{eq:FF}
	    \end{equation}
\item[(ii)] The state $\xb(t)$ of the system $\Sigma$ corresponding to the 
initial state $\xb(0)=x$ does not satisfy $\lim_{t\to\infty} \xb(t)\to 0$.  
\end{itemize}
\label{th:unstab}
\end{theorem}

For the proof of this theorem, we need the next lemma.
\begin{lemma}[\cite{Ebihara_2012}]
For given $a,b\in\bbR^n$ with $a\ne 0$, the following conditions are equivalent:
\begin{itemize}
 \item[(i)] $\He\{ab^T\}\succeq 0$.  
 \item[(ii)] There exists $\lambda\in\bbR^+$ such that $b=\lambda a$.  
\end{itemize}
\label{le:align}
\end{lemma}
\begin{proofof}{\rth{th:unstab}}
We first prove \eqref{eq:FF}.  
Since $\rank(H)=1$ holds, its full-rank factorization is given of the form
\[
   H=\begin{bmatrix}\hone\\\htwo\end{bmatrix}
 \begin{bmatrix}\hone\\\htwo\end{bmatrix}^T,\ 
 \hone\in\bbR^n,\ \htwo\in \bbR^m. 
\]
Then, from the second constraint in \eqref{eq:dLMI}, we have
\[
\begin{bmatrix} \htwo-(C\hone+D\htwo) \\ \htwo \end{bmatrix} 
\begin{bmatrix} \htwo-(C\hone+D\htwo) \\ \htwo \end{bmatrix}^T\ge 0.  
\]
The latter condition implies in particular that all entries have the same sign as $w_m$, the last entry of $w$. 
If $w_m \leq 0$, one can flip the sign of both $x$ and $w$ to obtain 
\begin{equation}
\begin{bmatrix} \htwo-(C\hone+D\htwo) \\ \htwo \end{bmatrix} \ge 0,
\label{eq:nonneg}
\end{equation}
yielding the rank one factorization \eqref{eq:FF}. 
\if{
Therefore we see 
\[
\begin{bmatrix} \htwo-(C\hone+D\htwo) \\ \htwo \end{bmatrix} \ge 0\ \mathrm{or} \ 
\begin{bmatrix} \htwo-(C\hone+D\htwo) \\ \htwo \end{bmatrix} \le 0,
\]
but this implies that (by changing the sign of $\hone$ and $\htwo$ if necessary)
we can always choose $\hone$ and $\htwo$ such that
\begin{equation}
\begin{bmatrix} \htwo-(C\hone+D\htwo) \\ \htwo \end{bmatrix} \ge 0.  
\label{eq:nonneg}
\end{equation}
This shows that \eqref{eq:FF} holds.  
}\fi

In addition, the third constraint in \eqref{eq:dLMI} reduces to
\begin{equation}
 (\htwo-(C\hone+D\htwo))\odot \htwo=0.   
\label{eq:Hadamard}
\end{equation}
From \eqref{eq:nonneg} and \eqref{eq:Hadamard} together with 
the semialgebraic set representation \eqref{eq:ReLU_alg}, 
we see
\begin{equation}
 \htwo=\Phi(C\hone+D\htwo).  
\label{eq:h_Phi}
\end{equation}
On the other hand, the first constraint in \eqref{eq:dLMI} reduces to
\begin{equation}
 \He\{(A\hone+B\htwo)\hone^T\}\succeq 0.  
\label{eq:h1h1}
\end{equation}
Here we first prove $\hone\ne 0$ by contradiction.  
To this end, suppose $\hone=0$. 
Then, from the underlying assumption that   
$H\ne 0$, we see $\htwo\ne 0$.  
Then, we see that \eqref{eq:Hadamard} reduces to 
$(\htwo-D\htwo)\odot \htwo=0$ and this implies
$\htwo^T\htwo=\htwo^TD\htwo$.  
However, this cannot hold for $\htwo\ne 0$
from our assumption that $\|D\|<1$.  
Then, by contradiction, we have $\hone\ne 0$.  

Since $\hone\ne 0$ as proved, we see from \rle{le:align} and 
\eqref{eq:h1h1} that there exists $\lambda\in\bbR_+$ such that
\begin{equation}
 A\hone+B\htwo=\lambda \hone.  
\label{eq:key} 
\end{equation}

We are now in the right position to complete the proof.  
To validate the assertion (ii), it suffices to show that 
the state trajectory $\xb(t)$ of the system $\Sigma$ given by \eqref{eq:Sigma}
with the initial condition $\xb(0)=\hone$
is given by $\xb(t)=e^{\lambda t}\hone$.  
We now do a reasoning very similar to the proof of Theorem \ref{th:gstab} with $\alpha = e^{\lambda t}$. 
From \eqref{eq:key}, \eqref{eq:h_Phi}, and \eqref{eq:co-linear}, we have
\begin{equation}
\begin{array}{@{}lcl}\displaystyle
 \frac{d}{dt} \xb(t)& = & e^{\lambda t}\lambda \hone\\
& = & e^{\lambda t} \left(A\hone+B\htwo\right)\\
& = & Ae^{\lambda t} \hone+Be^{\lambda t} \htwo\\
& = & A\xb(t)+Be^{\lambda t}\Phi(C\hone+D\htwo)\\     
& = & A\xb(t)+B\Phi(Ce^{\lambda t}\hone+De^{\lambda t}\htwo) \\    
& = & A\xb(t)+B\Phi(\zb(t)) \\    
\end{array}
\label{eq:delx}
\end{equation}
where
\[
\begin{array}{@{}lcl}\displaystyle
\zb(t)&=&Ce^{\lambda t}\hone+De^{\lambda t}\htwo\\
&=&C\xb(t)+De^{\lambda t}\Phi(C\hone+D\htwo)\\
&=&C\xb(t)+D\Phi(Ce^{\lambda t}\hone+De^{\lambda t}\htwo)\\
&=&C\xb(t)+D\Phi(\zb(t)).  
\end{array}
\]
It follows that $\zb(t)=(I-D \Phi)^{-1}(C\xb(t))$.  
By substituting this into \eqref{eq:delx}, we see that
$\xb(t)=e^{\lambda t}\hone$ surely satisfies \eqref{eq:Sigma}. 
Since $\lambda > 0$ and $\hone \neq 0$ the quantity $\xb(t)$ diverges. 
\end{proofof}
\begin{remark}
For robust stability analysis problems of 
linear dynamical systems affected by parametric uncertainties, 
it is well known that dual LMIs are useful 
in extracting worst-case parameters that destabilize the system,
see, e.g., 
\cite{Scherer_SIAM2005,Masubuchi_SCL2009,Ebihara_IEEE2009}.  
In the present result for the stability analysis of the nonlinear feedback system $\Sigma$, 
the dual LMI \eqref{eq:dLMI} works effectively 
in extracting a worst-case initial condition from which the
state trajectory does not converge to the origin.  
\end{remark}

\subsection{Numerical Examples}
\label{sec:num}

The systems considered in the examples in this section 
and \rsub{sub:num} are chosen to simulate the typical situation
of the stability analysis of RNNs and NN-driven control systems.  
If we employ a ReLU-RNN as a model of a complex dynamical system, 
we typically employ a large number of ReLUs to imitate the complex behavior.  
In addition, for a plant subject to unknown nonlinearities, 
we often employ a dynamical ReLU-NN as a controller 
again with a large number of ReLUs
to approximate the unknown nonlinearities.  
In such cases, the system of interest will be described 
by \eqref{eq:G} and \eqref{eq:Phi} where $m>n$.  

\subsubsection*{The Case $D\ne 0$ and feasible primal LMI}
\label{sub:Dne0Pfeas}

Let us consider the case where $n=2$, $m=5$ in \eqref{eq:G} and
\[
\arraycolsep=1mm
\scalebox{0.85}{$
\begin{array}{@{}l@{\hspace*{1mm}}l}
A=\begin{bmatrix*}[r]
 -1.36 & -0.53 \\
 -0.41 & -0.50 \\
\end{bmatrix*},\ 
&B=\begin{bmatrix*}[r]
 -0.82 & -0.79 &  0.57 &  0.21 & -0.14 \\
 -0.19 & -0.78 & -0.42 &  0.93 &  0.39 \\
  \end{bmatrix*} ,\vspace*{1mm}\\
C=\begin{bmatrix*}[r]
  0.52 & -0.63 \\
 -0.13 & -0.47 \\
  0.31 &  0.60 \\
 -0.78 & -0.02 \\
  0.87 &  0.54 \\
  \end{bmatrix*},\ 
&D=\begin{bmatrix*}[r]
 -0.21 & -0.10 &  0.39 & -0.15 &  0.26 \\
 -0.45 &  0.22 & -0.75 &  0.31 & -0.75 \\
 -0.93 & -0.88 & -0.74 &  0.45 & -0.73 \\
  0.35 & -0.37 & -0.82 &  0.06 & -0.80 \\
 -0.14 &  0.55 & -0.98 & -0.78 & -0.72 \\
  \end{bmatrix*}.  
\end{array}$}
\]
For this system, the primal LMI \eqref{eq:pLMI} turns out to feasible, 
thereby we can conclude that the system $\Sigma$ is stable.  
\rfig{fig:vector_field5.1} shows the vector field of system $\Sigma$, 
together with the state trajectories from the 
initial states $\xb(0)=[\ 3\ 3\ ]^T$ and $\xb(0)=[\ -3\ -3\ ]^T$.  
The state trajectories converge to the origin.  
In this example, the LMI \eqref{eq:stab} with the set of 
O'Shea-Zames-Falb multipliers employing
doubly hyperdominant matrices turns out to be infeasible.  
This shows the usefulness of the new set of multipliers
$\bPi_\NN$ proposed in \rpr{pr:PiCOP}.  

\subsubsection*{The Case $D=0$ and feasible dual LMI}
\label{sub:De0Dfeas}

Let us consider the case where $n=2$, $m=5$ in \eqref{eq:G} and
\[
\scalebox{0.85}{$
\begin{array}{@{}l@{\hspace*{1mm}}l}
A=\begin{bmatrix*}[r]
-1.35 & 0.7700 \\
-0.22 & 0.0400 \\
\end{bmatrix*},\ 
&B=\begin{bmatrix*}[r]
 -0.73 & -0.92 & 0.91 & 0.50 & -0.41 \\
 -0.90 & -0.55 & -0.13 & 0.05 & -0.09 \\
  \end{bmatrix*} ,\vspace*{1mm}\\
C=\begin{bmatrix*}[r]
  0.81 & -0.02 \\
 -0.11 & -0.35 \\
 -0.81 & 0.90 \\
  0.96 & -0.06 \\
  0.74 & 0.94 \\
  \end{bmatrix*},\ 
&D=0_{5,5} .  
\end{array}$}
\]
For this system, the dual LMI \eqref{eq:dLMI} turns out to be feasible, 
and the resulting dual solution $H$ is 
numerically verified to be $\rank(H)=1$.  
Therefore we can conclude that the system $\Sigma$ is not stable.  
The full-rank factorization of $H$ of the form \eqref{eq:FF} 
and $\lambda\in\bbR_+$ are given by
\[
 \hone=
 \begin{bmatrix*}[r]
   -0.6282\\
   -0.7780
 \end{bmatrix*},\ 
 \htwo=
 \begin{bmatrix*}[r]
   0.0000\\
   0.3414\\
   0.0000\\
   0.0000\\
   0.0000
\end{bmatrix*},\ \lambda=0.1037.  
\]
\rfig{fig:vector_field5.2} shows the vector field of system $\Sigma$, 
together with the state trajectories from 
initial states $\xb(0)=\hone$ and $\xb(0)=[\ 1\ 1\ ]^T$.  
The state trajectory from the initial state $\xb(0)=[\ 1\ 1\ ]^T$ converges 
to the origin.  
However, as proved in \rth{th:unstab}, the 
state trajectory from the initial state $\xb(0)=\hone$ does not converge
to the origin and in this case diverges along the ray $\xb(t) = e^{\lambda t} \hone$.  

\subsubsection*{The Case $D\ne 0$ and feasible dual LMI}
\label{sub:Dne0Dfeas}

\begin{figure}[htbp]
\centering
 \hspace*{-3mm}
 \includegraphics[scale=0.65]{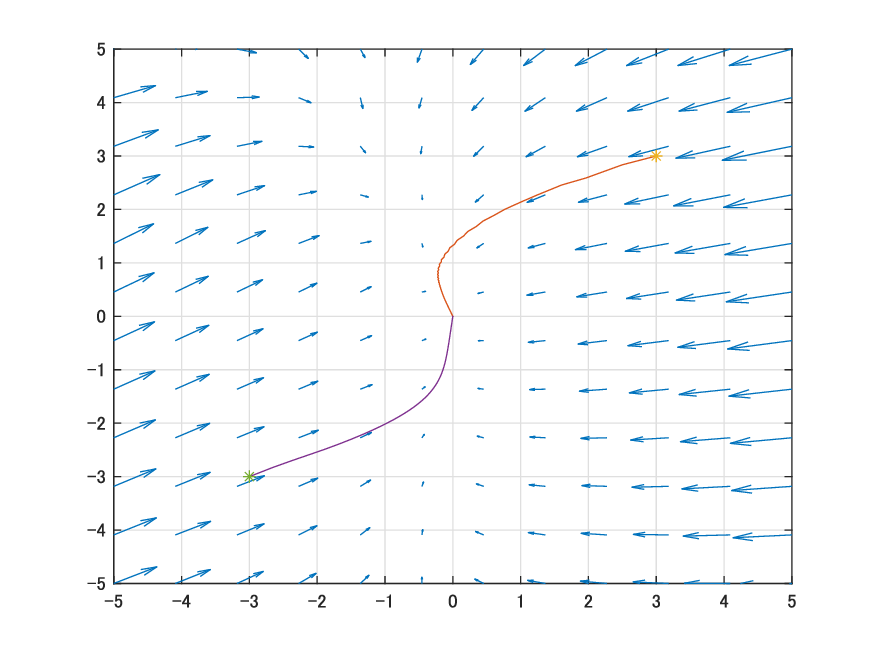}\vspace*{-3mm}
 \caption{State trajectories from initial conditions $\xb(0)=[\ 3\ 3\ ]^T$ and $\xb(0)=[\ -3\ -3\ ]^T$.  }
 \label{fig:vector_field5.1}
\centering
 \hspace*{-3mm}
 \includegraphics[scale=0.65]{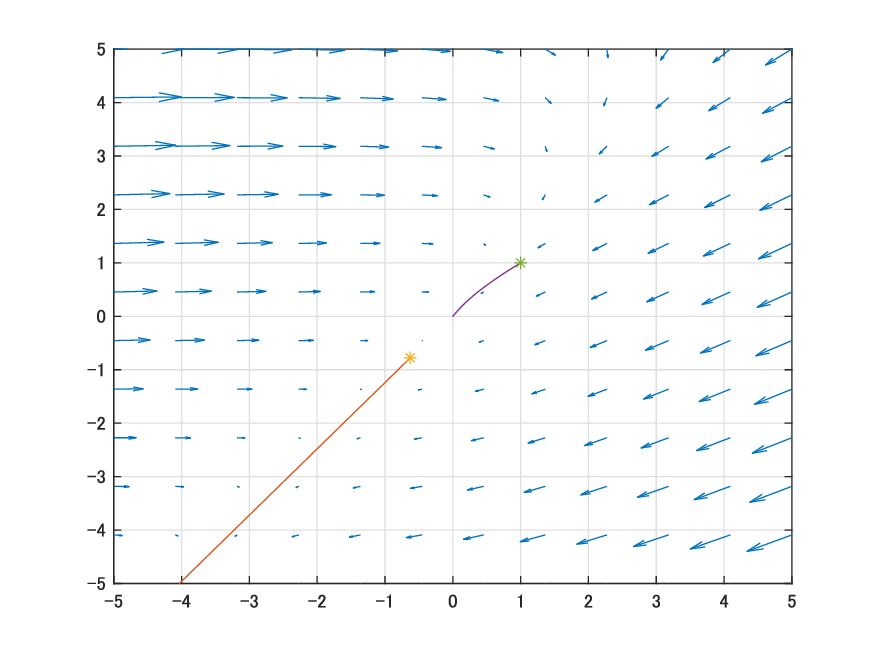}\vspace*{-3mm}
 \caption{State trajectories from initial conditions  $\xb(0)=\hone$ and $\xb(0)=[\ 1\ 1\ ]^T$.  }
 \label{fig:vector_field5.2}
\centering
 \hspace*{-3mm}
 \includegraphics[scale=0.65]{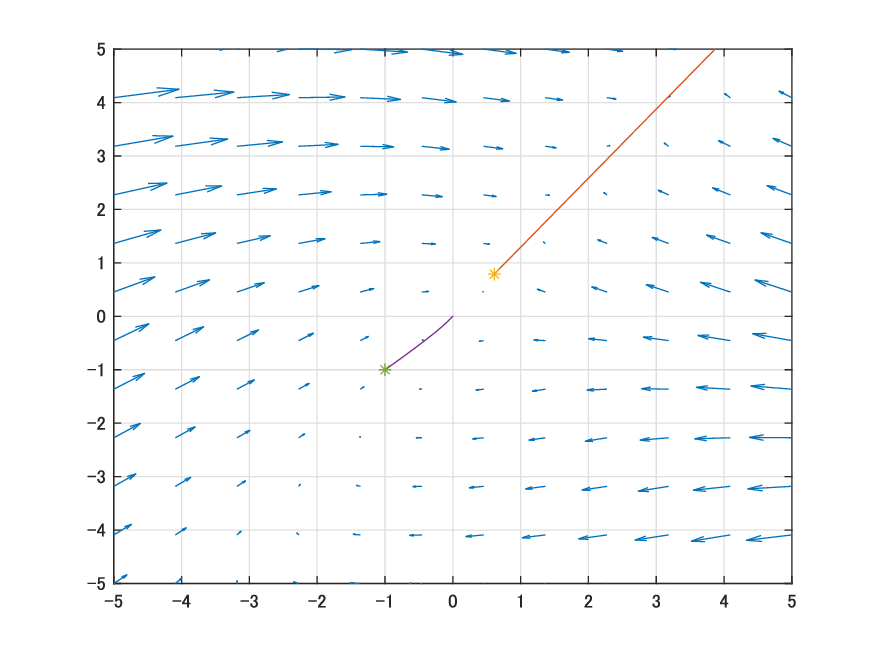}\vspace*{-3mm}
 \caption{State trajectories from initial conditions  $\xb(0)=\hone$ and $\xb(0)=[\ -1\ -1\ ]^T$.  }
 \label{fig:vector_field5.3}
\end{figure}

Let us consider the case where $n=2$, $m=5$ in \eqref{eq:G} and
\[
\scalebox{0.85}{$
\begin{array}{@{}l@{\hspace*{1mm}}l}
A=\begin{bmatrix*}[r]
 -1.17 & 0.86\\
 0.39 & -0.41\\
\end{bmatrix*},\ 
&B=\begin{bmatrix*}[r]
 -0.09 & 0.56 & -0.49 & 0.26 & 0.29\\
 -0.94 & -0.88 & 0.94 & 0.57 & 0.51\\
  \end{bmatrix*} ,\vspace*{1mm}\\
C=\begin{bmatrix*}[r]
 -0.40 & 0.44\\
 0.11 & -0.90\\
 -0.50 & -0.07\\
 -0.88 & 0.61\\
 0.18 & 0.35\\
  \end{bmatrix*},\ 
&D=\begin{bmatrix*}[r]
 -0.38 & 0.15 & -0.14 & -0.27 & -0.37\\
 0.14 & 0.19 & -0.19 & 0.16 & -0.17\\
 0.34 & 0.36 & 0.05 & 0.06 & -0.23\\
 -0.37 & -0.17 & 0.16 & -0.17 & -0.25\\
 0.53 & -0.40 & -0.33 & -0.08 & -0.32\\
  \end{bmatrix*}.  
\end{array}$}
\]
For this system, the dual LMI \eqref{eq:dLMI} turns out to be feasible, 
and the resulting dual solution $H$ is 
numerically verified to be $\rank(H)=1$.  
Therefore we can conclude that the system $\Sigma$ is not stable.  
The full-rank factorization of $H$ and $\lambda\in\bbR_+$ are given by
\[
 \hone=
 \begin{bmatrix*}[r]
    0.6119\\
    0.7909\\
 \end{bmatrix*},\ 
 \htwo=
 \begin{bmatrix*}[r]
   0.0000\\
   0.0000\\
   0.0000\\
   0.0000\\
   0.2932\\
\end{bmatrix*},\ \lambda=0.0807.  
\]
\rfig{fig:vector_field5.3} shows the vector field of system $\Sigma$, 
together with the state trajectories from 
initial states $\xb(0)=\hone$ and $\xb(0)=[\ -1\ -1\ ]^T$.  
The state trajectory from the initial state $\xb(0)=[\ -1\ -1\ ]^T$ converges 
to the origin.  
However, as proved in \rth{th:unstab}, the 
state trajectory from the initial state $\xb(0)=\hone$ does not converge
to the origin (and in this case diverges).  

We finally note that, from the comparison between 
Figs. \ref{fig:vector_field5.1}, \ref{fig:vector_field5.3} and 
\rfig{fig:vector_field5.2}, 
we can recognize that the appearance of the vector fields
with $D\ne 0$ becomes more complicated 
than that of $D=0$.  

%

\section{Higher order LMI relaxations}
\label{sec:higher}
\subsection{Dual LMI interpretation}
\label{sec:interpret}
The dual LMI \eqref{eq:dLMI} and the proof of \rth{th:unstab}
together with the semialgebraic set representation \eqref{eq:ReLU_alg} 
motivate us to consider the following 
algebraic inequality feasibility problem:

\noindent
Find nonzero $\hone\in\bbR^n,\ \htwo\in\bbR^m,\ \lambda \in \bbR_+$ such that
\begin{equation}
\begin{array}{@{}l}
A\hone+B\htwo=\lambda \hone,\\ 
\begin{bmatrix} \htwo-(C\hone+D\htwo) \\ \htwo \end{bmatrix}\ge 0,\\ 
 (\htwo-(C\hone+D\htwo))\odot \htwo=0.   
\end{array}
\label{eq:alge}
\end{equation}

For this problem we readily obtain the next result.  
\begin{theorem}
For given $\hone\in\bbR^n,\ \htwo\in\bbR^{m}$ and $\lambda \in \bbR_+$, 
the following conditions are equivalent.
\begin{itemize}
 \item[(i)] The condition \eqref{eq:alge} holds.  
 \item[(ii)] The solution of the system $\Sigma$ given by \eqref{eq:Sigma} with $\xb(0)=\hone$
	     is given by $\xb(t)=e^{\lambda t} \hone$.  
\end{itemize}
\label{th:equiv_sol}
\end{theorem}
\begin{proofof}{\rth{th:equiv_sol}}
We have showed (i) $\Ra$ (ii) in the proof of \rth{th:unstab}.  
Therefore it suffices to prove (ii) $\Ra$ (i).  
To this end, suppose (ii) holds.  Then, we have from \eqref{eq:Sigma} that
\[
\begin{array}{@{}lcl}
e^{\lambda t}  \lambda \hone&=&Ae^{\lambda t}\hone+B\Phi\circ(I-D\Phi)^{-1}(Ce^{\lambda t}\hone)\\
&=&e^{\lambda t} \left(A\hone+B\Phi\circ(I-D\Phi)^{-1}(C\hone)\right)\\
 \end{array}
\]
where we used \eqref{eq:co-linear} and \eqref{eq:co-linear-invmap}.  
Then, if we define 
\begin{equation}
\htwo:=\Phi\circ(I-D\Phi)^{-1}(C\hone)\in\bbR_+^m,  
\label{eq:h2}
\end{equation}
we have
\begin{equation}
 A\hone+B\htwo=\lambda \hone.    
\label{eq:proof1} 
\end{equation}
Here, if we let $\zeta:=(I-D\Phi)^{-1}(C\hone)$, 
we see that
\[
\scalebox{1.0}{$
\begin{array}{@{}l}
(I-D\Phi)(C\hone+D\htwo)\\
=(I-D\Phi)(C\hone+D\Phi\circ(I-D\Phi)^{-1}(C\hone))\\
=(I-D\Phi)((I-D\Phi)(\zeta)+D\Phi(\zeta))\\
=(I-D\Phi)(\zeta)\\
=C\hone.  
\end{array}$}
\]
It follows that we have in fact
\[
(I-D\Phi)^{-1}(C\hone)=C\hone+D\htwo.  
\]
This, together with \eqref{eq:h2}, leads us to
\[
 \htwo=\Phi(C\hone+D\htwo).    
\]
From \eqref{eq:ReLU_alg}, we therefore obtain
\begin{equation}
\begin{array}{@{}l}
\begin{bmatrix} \htwo-(C\hone+D\htwo) \\ \htwo \end{bmatrix}\ge 0,\\ 
 (\htwo-(C\hone+D\htwo))\odot \htwo=0.   
\end{array}
\label{eq:proof2}
\end{equation}
From \eqref{eq:proof1} and \eqref{eq:proof2}, we can conclude that 
(i) holds.  This completes the proof.  
\end{proofof}

The algebraic inequality feasibility problem \eqref{eq:alge}
is nonconvex and hence numerically intractable.  
In fact, if we follow \cite{Ebihara_IEEE2009}, 
the dual LMI \eqref{eq:dLMI} can be interpreted as 
the first order LMI relaxation of this nonconvex feasibility problem.  
To see this, for the variables 
$\hone\in\bbR^n,\ \htwo\in\bbR^{m}$ and $\lambda \in \bbR_+$, 
let us define the matrix variable
\begin{equation}
\tilde{H}_1=
\begin{bmatrix} H_{11} & H_{12} \\ H_{12}^T & H_{22} \end{bmatrix}
=
\begin{bmatrix} \hone \\ \htwo \end{bmatrix}
\begin{bmatrix} \hone \\ \htwo \end{bmatrix}^T\in \bbS_+^{n+m}  
\label{eq:rank-one} 
\end{equation}
where the variable $\lambda \in \bbR_+$ is hidden but shows up via \rle{le:align}.  
Then, from \rle{le:align}, we see that
the first equality condition in \eqref{eq:alge} is equivalent to
the first inequality condition in \eqref{eq:dLMI}, 
and the second inequality and 
the third equality conditions in \eqref{eq:alge} are equivalent
to the second inequality and the third equality conditions in \eqref{eq:dLMI}, respectively.  
By removing the nonconvex rank-one constraint \eqref{eq:rank-one}, 
we then obtain the convex LMI relaxation problem \eqref{eq:dLMI}.  
With this interpretation in mind, we now present higher order relaxations
in the next subsection.  

\subsection{Higher order LMI relaxations}
\label{sec:subhigh}
For general polynomial optimization problems, 
a hierarchy of LMI relaxations with asymptotic exactness 
is proposed in \cite{lasserre2001global}.   
Here we similarly derive the so-called \textit{second order} LMI relaxation for the problem \eqref{eq:alge}.  
This relaxation is similar to the one from 
\cite{Ebihara_IEEE2009}, and is given as follows:

\noindent
{\bf Second Order LMI Relaxation}
\begin{equation}
\scalebox{0.95}{$
\begin{array}{@{}l}
 \mathrm{Find}\ \tilde{H}_2=
\begin{bmatrix} 
H_{11} & H_{12} & H_{13} & H_{14} \\ 
\ast      & H_{22} & H_{14}^T & H_{24} \\ 
\ast      & \ast      & H_{33} & H_{34} \\ 
\ast      & \ast      & \ast      & H_{44} \\ 
\end{bmatrix}
\in\bbS_+^{2(n+m)}\setminus \{0\}\\ 
\mathrm{such\ that}\\
\begin{bmatrix} 
H_{13} & H_{14} \\ 
H_{14}^T & H_{24} 
\end{bmatrix}\in \bbS_+^{n+m},\\
AH_{11}+BH_{12}^T=H_{13},\\ 
AH_{12}+BH_{22}=H_{14},\\ 
AH_{13}+BH_{14}^T=H_{33},\\ 
AH_{14}+BH_{24}=H_{34},\\ 
\He\{AH_{33}+BH_{34}^T\}\succeq 0,\\ 
\begin{bmatrix} -C & -D+I_{m} \\ 0 & I_{m} \end{bmatrix} 
\begin{bmatrix} H_{11} & H_{12} \\ \ast & H_{22}\end{bmatrix} 
\begin{bmatrix} -C & -D+I_{m} \\ 0 & I_{m} \end{bmatrix}^T \ge 0,\\ 
\begin{bmatrix} -C & -D+I_{m} \\ 0 & I_{m} \end{bmatrix} 
\begin{bmatrix} H_{13} & H_{14} \\ \ast & H_{24}\end{bmatrix} 
\begin{bmatrix} -C & -D+I_{m} \\ 0 & I_{m} \end{bmatrix}^T \ge 0,\\ 
\begin{bmatrix} -C & -D+I_{m} \\ 0 & I_{m} \end{bmatrix} 
\begin{bmatrix} H_{33} & H_{34} \\ \ast & H_{44}\end{bmatrix} 
\begin{bmatrix} -C & -D+I_{m} \\ 0 & I_{m} \end{bmatrix}^T \ge 0,\\ 
\diag(-CH_{12}+(-D+I_m)H_{22}) =0,\\  
\diag(-CH_{14}+(-D+I_m)H_{24}) =0,\\  
\diag(-CH_{34}+(-D+I_m)H_{44}) =0.  
\end{array}$}
\label{eq:dLMI2}
\end{equation}
Similarly, by defining 
$\clIu:=[\ I_n \ 0_{n,m}\ ]\in\bbR^{n\times (n+m)}$ and 
$\clIl:=[\ 0_{m,n}\ I_m\ ]\in\bbR^{m\times (n+m)}$, 
the $N$-th order LMI relaxation for the problem \eqref{eq:alge} can be given as follows: 

\noindent
{\bf $N$-th Order LMI Relaxation}
\begin{equation}
\scalebox{0.95}{$
\begin{array}{@{}l}
 \mathrm{Find}\ \tilde{H}_N=
\begin{bmatrix} 
\clH_{0} & \clH_{1} & \cdots  & \clH_{N-1} \\ 
\clH_{1} & \iddots  & \iddots & \vdots \\ 
\vdots    & \iddots  & \iddots & \clH_{2N-3} \\
\clH_{N-1}  & \cdots & \clH_{2N-3} & \clH_{2(N-1)}
\end{bmatrix}
\in\bbS_+^{N(n+m)}\setminus \{0\}\hspace*{-20mm}\\ 
\mathrm{such\ that}\\
\clH_{i}\in \bbS_+^{n+m}\ (i=0,\cdots,2(N-1)),\\
\begin{bmatrix} A & B \end{bmatrix} \clH_{i} = \clIu \clH_{i+1} \ (i=0,\cdots,2N-3),\\
\He\{\begin{bmatrix} A & B \end{bmatrix} \clH_{2(N-1)}\clIu^T\}\succeq 0,\\ 
\begin{bmatrix} -C & -D+I_{m} \\ 0 & I_{m} \end{bmatrix} 
\clH_i
\begin{bmatrix} -C & -D+I_{m} \\ 0 & I_{m} \end{bmatrix}^T \ge 0\\ 
(i=0,\cdots,2(N-1)),\\
\diag(\begin{bmatrix} -C & -D+I_m \end{bmatrix}\clH_i \clI_l^T)=0 \ (i=0,\cdots,2(N-1)).  \hspace*{-20mm}
\end{array}$}
\label{eq:dLMIN}
\end{equation}
We now make several observations about 
the $N$-th order LMI relaxation \eqref{eq:dLMIN}. 
\begin{enumerate}[leftmargin=*]
\item  The matrix variable $\tilde{H}_N \in\bbS^{N(n+m)}$ has
       a block-Hankel matrix structure.  
       A matrix variable of this form has already been employed in \cite{Ebihara_IEEE2009} 
       to deal with robust stability analysis problems of linear dynamical systems
       affected by a single uncertain parameter.  
       In relation to the variables $\hone\in\bbR^n$, $\htwo\in\bbR_+^m$, and 
       $\lambda \in \bbR_+$ of the original algebraic inequality feasibility problem 
       \eqref{eq:alge}, 
       the matrix variable $\tilde{H}_N\in\bbS_+^{N(n+m)}$ corresponds to the relaxation
       of the rank-one matrix variable
       \begin{equation}
       \begin{bmatrix} \hone \\ \htwo\\ \vdots \\ \lambda^{N-1} \hone \\ \lambda^{N-1} \htwo  \end{bmatrix}
       \begin{bmatrix} \hone \\ \htwo\\ \vdots \\ \lambda^{N-1} \hone \\ \lambda^{N-1} \htwo  \end{bmatrix}^T\in \bbS_+^{N(n+m)}.  
       \label{eq:Nthval}
       \end{equation}
       Therefore the feasibility of \eqref{eq:dLMIN} is a necessary condition
       for the feasibility of \eqref{eq:alge}. 
\item  Suppose $\rank(\tilde{H}_N)=1$ in the $N$-th order relaxation \eqref{eq:dLMIN}.  
       Then,  the system $\Sigma$ is unstable.  
       Moreover, the full-rank factorization 
       of $\tilde{H}_N$ is given of the form \eqref{eq:Nthval} with $\lambda\ge 0$, 
       and the same assertion with (ii) of \rth{th:unstab} holds.  
\item For $N_1\le N_2$, suppose \eqref{eq:dLMIN} is feasible for $N_2$.  
       Then, from the structure of the constraints in \eqref{eq:dLMIN},
       it is very clear that \eqref{eq:dLMIN} is feasible for $N_1$.  
       Conversely, from \eqref{eq:Nthval}, it is also true that, 
       if $\tilde{H}_{N_1} \in\bbS_+^{N_1(n+m)}$ with $\rank(\tilde{H}_{N_1})=1$
       satisfies \eqref{eq:dLMIN}, then there exists 
       $\tilde{H}_{N_2} \in \bbS_+^{N_2(n+m)}$ with $\rank(\tilde{H}_{N_2})=1$
       that satisfies \eqref{eq:dLMIN}.  
       Therefore, by increasing the order $N$, 
       we can appropriately restrict the dual solution so that the satisfaction of 
       the rank-one condition is more likely to happen.  
 \item As stated above, it is promising to carry out higher order relaxations to obtain
       instability certificate for the system $\Sigma$.  
       We provide later on in Section \ref{sec:setting} another related hierarchy of primal-dual LMI relaxation with asymptotic convergence guarantees with respect to the original feasibility problem \eqref{eq:alge}. 
       %
 \item If we follow the standard and well-established 
       LMI relaxation methodology for general polynomial optimization problems
       \cite{lasserre2001global}, we employ the (relaxed one of the) following so-called \textit{moment} matrix
       as the variable in the first order LMI relaxation of \eqref{eq:alge}:
       \[
       \begin{bmatrix} 1 \\ \lambda\\ \hone \\ \htwo \end{bmatrix}
       \begin{bmatrix} 1 \\ \lambda\\ \hone \\ \htwo \end{bmatrix}^T
       \in \bbS_+^{n+m+2}.    
       \]
       This moment matrix variable includes
       all the monomials with respect to the variables
       $\lambda \in \bbR,\ \hone\in\bbR^n$, and $\htwo\in\bbR^m$
       of degrees up to two.  
       In stark contrast, in the present first order LMI relaxation \eqref{eq:dLMI}, 
       the dual variable only includes the monomials of degree two
       with respect to the variables $\hone\in\bbR^n$ and $\htwo\in\bbR^m$ 
       as shown in \eqref{eq:rank-one}. In this sense, it is sparse.  
       Similarly, in the $N$-th order moment relaxation, 
       the size of the moment matrix variable employed
       is $\binom{n+m+1+N}{N}$ and it includes 
       all the monomials with respect to the variables
       $\lambda \in \bbR,\ \hone\in\bbR^n,\ \htwo\in\bbR^m$
       of degrees up to $N$.  
       However, in the above $N$-th order relaxation, 
       the size of the matrix variable is $N(n+m)$ as shown in \eqref{eq:Nthval} and 
       it is sparse.  
       Convergence proofs have been previously obtained 
       while exploiting either \emph{correlative} or \emph{term} sparsity; 
       see \cite{Magron_2023} for a recent monograph on the topic. 
       In the present context, the monomial choice in our sparse hierarchy of 
       dual LMIs does not directly fit within previously studied frameworks, 
       which motivates the  additional investigation from Section \ref{sec:setting}. 

%

\item We finally note that, 
to make the rank-one condition more likely to be satisfied, 
we apply the following well-known heuristic, mentioned, e.g., in \cite{Henrion_2005b}: 
\[
 \mathrm{min}\ \trace(\tilde{H}_N)\ \mathrm{s.t.}\ \eqref{eq:dLMIN}\ \mathrm{and}\ 
\trace(\clIu\clH_0\clIu^T)=1.  
\]
With this trace minimization, the rank of $\tilde{H}_N$ tends to be reduced
and thus we can expect that the suggested rank-one condition 
becomes more likely to be satisfied.  
This heuristic has been also employed 
for the numerical examples presented in Section \ref{sub:De0Dfeas}.  
\end{enumerate}
\subsection{Numerical examples}
\label{sub:num}

\begin{figure}[t]
\centering
 \hspace*{-3mm}
 \includegraphics[scale=0.65]{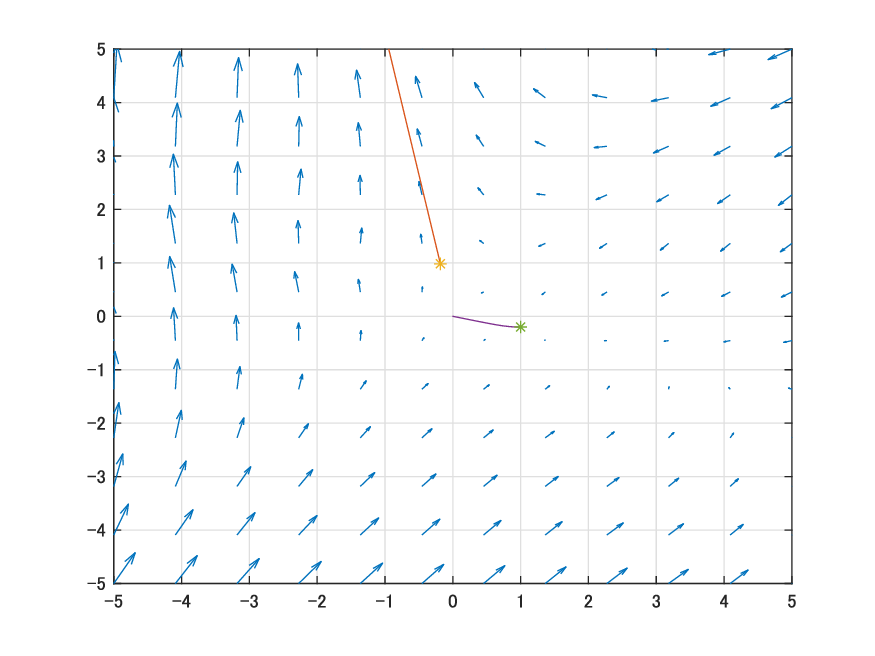}\vspace*{-3mm}
 \caption{State trajectories from the initial conditions  $\xb(0)=\hone$ and $\xb(0)=[\ 1\ -0.2\ ]^T$.  }
 \label{fig:vector_field6.3}
\end{figure}

Let us consider the case where $n=2$, $m=5$ in \eqref{eq:G} and
\[
\scalebox{0.85}{$
\begin{array}{@{}l@{\hspace*{1mm}}l}
A=\begin{bmatrix*}[r]
 -0.21 &  0.00\\
 -0.14 & -0.23\\
\end{bmatrix*},\ 
&B=\begin{bmatrix*}[r]
  0.14 & -0.67 & -0.37 &  0.46 &  0.11\\
  0.15 &  0.92 &  0.17 & -0.01 &  0.77\\
  \end{bmatrix*} ,\vspace*{1mm}\\
C=\begin{bmatrix*}[r]
  0.65 & -0.61\\
 -0.26 &  0.50\\
 -0.72 & -0.46\\
  0.11 & -0.65\\
 -0.60 &  0.61\\
  \end{bmatrix*},\ 
&D=\begin{bmatrix*}[r]
 -0.18 &  0.24 & -0.15 & -0.32 & -0.36\\
  0.19 & -0.38 &  0.21 &  0.03 & -0.28\\
  0.38 &  0.04 &  0.14 & -0.40 & -0.34\\
 -0.23 & -0.05 &  0.27 &  0.14 &  0.28\\
  0.20 &  0.08 & -0.32 &  0.18 & -0.34\\
  \end{bmatrix*}.  
\end{array}$}
\]
For this system, the primal LMI \eqref{eq:pLMI} turns out to be infeasible, 
and thus we cannot conclude anything about the stability of $\Sigma$ at this stage.  
Therefore we apply the hierarchy of LMI relaxations shown in \eqref{eq:dLMIN}.  
The first order LMI relaxation of course turns out to be feasible, 
but the resulting dual variable $\clH$ does not satisfy $\rank(\clH)=1$.  
Similarly for the second order LMI relaxation.  
However, the third order LMI relaxation turns out to be feasible, 
and in particular the resulting dual variable $\clH$ does satisfy $\rank(\clH)=1$.  
The full-rank factorization of $\clH$ of the form \eqref{eq:Nthval} 
and $\lambda\in\bbR_+$ are given by
\[
 \hone=
 \begin{bmatrix*}[r]
   -0.1831\\
    0.9831\\
 \end{bmatrix*},\ 
 \htwo=
 \begin{bmatrix*}[r]
   0.0000\\
   0.2799\\
   0.0000\\
   0.0000\\
   0.5462\\
\end{bmatrix*},\ \lambda=  0.4858.  
\]
\rfig{fig:vector_field6.3} shows the vector field of system $\Sigma$, 
together with the state trajectories from 
initial states $\xb(0)=\hone$ and $\xb(0)=[\ 1\ -0.2\ ]^T$.  
The state trajectory from the initial state $\xb(0)=[\ 1\ -0.2\ ]^T$ converges 
to the origin.  
However, as shown at the end of Section \ref{sec:subhigh}, the 
state trajectory from the initial state $\xb(0)=\hone$ does not converge
to the origin, and actually diverges in this case.  

\section{Complete LMI hierarchies for instability analysis}
\label{sec:setting}
Let us consider the system $\Sigma$ as in \eqref{eq:G} with $A \in \R^{n \times n}, B \in \R^{n \times m}, C \in \R^{m \times n}$ and $D \in \R^{m \times m}$. 
Based on the considerations from Section \ref{sec:interpret}, the instability of $\Sigma$ can be certified by considering the feasibility problem stated in \eqref{eq:alge}. \\
In Section \ref{sec:polsys}, we consider a specific related system of polynomial (in)equalities. 
We show in Section \ref{sec:exploit} how to exploit the structure of equality constraints. 
Then we derive in Section \ref{sec:complete} a \textit{complete} hierarchy of LMIs that allows one to prove that this system is feasible. 
Here ``complete'' means that the hierarchy shall always converge towards a limit where one can assert the desired instability statement. 
In Section \ref{sec:micnon}, we relate more precisely this hierarchy with the one provided in Section \ref{sec:higher}. 
\subsection{Preliminaries on polynomial systems}
\label{sec:polsys}
%
We consider the following feasibility problem: \\
Find $x \in \R^n, w \in \R^m, \lambda \in \R$ such that
\begin{align}
A x + B w - \lambda x = 0, \label{eqlinear}\\
w \odot (C x + D w - w) = 0, \label{eqrelu}\\
\sum_{i=1}^n x_i^2 + \sum_{j=1}^m w_j^2 -1 = 0, \label{equnit}\\
\lambda \geq 0, \quad g(x, w) \geq 0, \quad \forall g \in G \label{ineqquadr},
\end{align}
where $G$ is the set of polynomials being quadratic in $(x, w)$ given by 
\begin{align*}
G = &\left\lbrace w_i w_j, w_j w_i - w_j (C x + D w)_i,  \right.\\
& \left. [w_i - (C x + D w)_i] [w_j - (C x + D w)_j] \mid 1\leq i,j \leq m \right\rbrace. 
\end{align*}
Let $\R[x, w,\lambda]$ be the vector space of polynomials in variables $(x, w,\lambda)$. 
Let $E := (e_i)_{1\leq i \leq n+m+1}$ be the set of $n+m+1$ polynomials involved in \eqref{eqlinear},  \eqref{eqrelu} and \eqref{equnit}, namely $e_i := (A x + B w - \lambda x)_i$, $i =1,\dots,n$,  $e_{n+i} := w_i (C x + D w)_i - w_i^2$, $i =1,\dots,m$ and $e_{n+m+1} := \sum_{i=1}^n x_i^2 + \sum_{j=1}^m w_j^2 -1$. 
Let us consider the \textit{ideal} $I$ generated by $E$: 
\begin{align}
\label{eq:ideal}
I := \left\lbrace \sum_{i=1}^{n+m+1} p_i e_{i} : p_i \in \R[x,w,\lambda] \right\rbrace. 
\end{align}
The quotient space $\R[x, w,\lambda]/I$ is the set of elements of the form $f + I = \{f + h \mid h \in I\}$, for $f \in \R[x, w,\lambda]$. 
The \textit{complex variety} associated to I is 
\begin{align}
\label{eq:variety}
V_{\C}(I) := \{(x, w,\lambda) \in \C^{n+m+1} \mid h(x, w,\lambda) = 0  \quad \forall h \in I \}.
\end{align}
When $V_{\C}(I)$ is finite then the ideal $I$ is said to be \emph{zero-dimensional}. 
Moreover the set 
\begin{align}
\label{eq:radical}
\sqrt{I} := \{h \in \R[x, w,\lambda] \mid h^k \in I \text{ for some integer } k \geq 1 \}
\end{align}
is the \emph{radical} of $I$. 
The ideal $I$ is said to be a radical ideal when $I = \sqrt{I}$. 
When $I$ is radical, $\# V_{\C}(I) = \dim \R[x, w,\lambda]/I $ and each root in $V_{\C}(I)$ has single multiplicity. 

Let us denote by `$<_{\lex}$' the lexicographical order on $\N^{n+m+1}$ such that $\lambda <_{\lex} x_1<_{\lex} \cdots <_{\lex} x_n <_{\lex} w_1<_{\lex} \cdots <_{\lex} w_n$. 
Given a polynomial $p \in \R[x, w,\lambda]$, one can write 
\[
p = \sum_{\alpha} p_\alpha x_1^{\alpha_1} \dots x_n^{\alpha_n} w_1^{\alpha_{n+1}} \dots w_m^{\alpha_{n+m}} \lambda^{\alpha_{n+m+1}}
\] 
and the \emph{leading term} $\LT(p)$ is defined as the maximum 
\[
p_{\alpha} x_1^{\alpha_1} \dots x_n^{\alpha_n} w_1^{\alpha_{n+1}} \dots w_m^{\alpha_{n+m}} \lambda^{\alpha_{n+m+1}}
\] 
with respect to $<_{\lex}$. 
The \emph{leading term ideal} of $I$ is $\LT(I) = \{\LT(h) \mid h \in I \}$. \\
The set of \emph{standard monomials} $\B$ with respect to $<_{\lex}$ is the set of monomials that do not belong to $LT(I)$, i.e., the set of monomials that cannot be divided by any term of the form $\LT(h)$ for $h \in I$. 
A finite subset $\Gb \subseteq I$ is called a \emph{Gr\"{o}bner basis} of $I$ if $\LT(I)=\LT(\Gb)$; that is, if the leading term  of every nonzero polynomial in $I$ is divisible by the leading term of some polynomial in $\Gb$. 
The existence of a Gr\"{o}bner basis is always ensured. A Gr\"{o}bner basis is called \emph{reduced} under the two following conditions: (1) the leading term coefficient of every element of $\Gb$ is 1 and (2) 
for all $g \in \Gb$, there are no monomials of $g$ lying in $\langle \LT(\Gb) \backslash \{g\}  \rangle $. Every ideal $I$ has a unique reduced Gr\"{o}bner basis. 
We refer the interested reader to \cite{cox2013} for more details.

Assume that $J \subseteq \R[x,y,\lambda]$ is a zero-dimensional and radical ideal and that $\Gb$ is the reduced Gr\"{o}bner basis of $J$ with respect to $<_{\lex}$. 
One says that $J$ is in \textit{shape position} if $\Gb$ has the following form:
\begin{equation}
\label{f:sysshape}
\Gb=\{u_0,x_1-u_1,\dots,x_n-u_n,w_1-v_1,\dots, w_m - v_m\},
\end{equation}
where $u_i,v_j$ are polynomials in $\R[\lambda]$ and $\deg u_0=\# V_{\C}(J)$.
We remind the so-called \emph{Shape Lemma} providing a criteria for an ideal to be in shape position. 

\begin{lemma}[Shape Lemma, \cite{gianni89}]
\label{lm:shape}
Let $J$ be a zero-dimensional and radical ideal and $<_{\lex}$ be a lexicographic monomial order as above. If $V_{\C}(J)$ is the union of points in $\C^{n+m+1}$ with distinct $\lambda$-coordinates, then $I$ is in shape position as in \eqref{f:sysshape}.
\end{lemma}
Eventually, we define the suitable notion of \textit{genericity} that shall ensure that the ideal $I$ is zero-dimensional and radical. 
For any $\cJ \subseteq \{1,\dots,m\}$, let us denote by $D_{\cJ \times \cJ}$ the principal sub-matrix of $D$ obtained by selecting only rows and columns indexed by $\cJ$. 
The number of such principal sub-matrices is equal to $2^m$. 
For any principal sub-matrix $D_{\cJ \times\cJ}$ of $D$, observe that $\|D_{\cJ \times \cJ}\| \leq \|D\|$ and thus that $I - D_{\cJ \times \cJ}$ is invertible. 
Let $F_{\cJ} \in \R^{m \times m}$ be the matrix such that its principal sub-matrix obtained by selecting only rows and columns indexed by $\cJ$ is $(I - D_{\cJ \times \cJ})^{-1}$, and with zero entries otherwise. 


\begin{definition}
\label{def:generic}
A matrix tuple $(A, B, C, D) \in \R^{n \times n} \times \R^{n \times m} \times \R^{m \times n} \times \R^{m \times m}$ with $\|D\| < 1$ is called $\textit{generic}$ if the following hold.
\begin{enumerate} 
\item For all $\cJ \subseteq \{1,\dots,m\}$, the matrix $I + F_{\cJ} C$ is invertible and the matrix $(I + F_{\cJ})(A+B F_{\cJ} C) (I + F_{\cJ} C)^{-1}$ has $n$ distinct nonzero eigenvalues, each of them being associated to two unit eigenvectors. Let us note $\Lambda_{\cJ}$ and $\mathcal{V}_{\cJ}$ the corresponding sets of eigenvalues and eigenvectors; 
\item The set $\bigcup_{\cJ \subseteq \{1,\dots,m\}} \Lambda_{\cJ}$ has $N = n  2^m$ distinct elements;
\item The set $ \bigcup_{\cJ \subseteq \{1,\dots,m\}} \{(I + F_{\cJ} C) v :  v \in \mathcal{V}_{\cJ} \}$ has $2 N $ distinct elements with nonzero first coordinates. 
\end{enumerate}
\end{definition}
\subsection{Algebraic structure exploitation}
\label{sec:exploit}
\begin{proposition}
\label{prop:zerodim}
Let us assume that the matrix tuple $(A, B, C, D) \in \R^{n \times n} \times \R^{n \times m} \times \R^{m \times n} \times \R^{m \times m}$, with $\|D\| < 1$, is generic in the sense of Definition \ref{def:generic}. 
Then, the ideal $I$ defined in \eqref{eq:ideal} is zero-dimensional, radical, and has a reduced Gröbner basis $\Gb$ with respect to $<_{\lex}$ of the following form: 
\begin{equation}
\label{eq:standard}
\Gb=\{u_0, x_1^2 - u_1, x_2-u_2 x_1,\dots,x_n-u_n x_1,w_1-v_1 x_1,w_2 - v_2 x_1, \dots, w_m - v_m x_1\},
\end{equation}
where $u_i,v_j$ are polynomials in $\R[\lambda]$ and $\deg u_0=\# V_{\C}(I)$.
\end{proposition}
\begin{proofof}{\rpr{prop:zerodim}}
We will show that the system of polynomial equations $e_1= \dots =  e_{n+m+1} = 0$, that involves $n+m+1$ variables $x, w,\lambda$, has finitely many complex solutions (possibly no real ones). 
Each equality constraint from \eqref{eqrelu} is equivalent to $w_i=0$ or $(C x + D w)_i - w_i = 0$. 
Let $\cJ \subseteq \{1,\dots,m\}$ be the set of indices corresponding to nonzero entries of $w$. 
\if{
For instance, if $w=0$ then the solutions of the system correspond to the unit eigenvectors of $A$ and associated eigenvalues. 
In this case one has $\cJ = \emptyset$, $F_{\cJ} = 0$. By the genericity assumption, the matrix $A$ has $n$ distinct eigenvalues, the set of unit eigenvectors is finite (2 for each eigenvalue, generically). 
The exact same reasoning holds if we assume that $w\neq 0$. 
}\fi
Then \eqref{eqrelu} is equivalent to $ w = F_{\cJ} C x$.  
We give two simple illustrative examples. 
If $w=0$, then $\cJ = \emptyset$, $F_{\cJ} = 0$ and the solutions of the system correspond to the unit eigenvectors of $A$ and associated eigenvalues. 
If $w_1 = 0$ and $w_i \neq 0$ for all $i>1$ then $F_{\cJ} = \begin{pmatrix}
0 & 0 \\
0 & (I - D_{m-1})^{-1}
\end{pmatrix}$ where $D_{m-1}$ is the matrix obtained by selecting entries from the $(m-1)$ last rows and columns of $D$. \\
The initial system of equations is equivalent to $(A + B F_{\cJ} C) x = \lambda x$, $ \|I + F_{\cJ} C x \|^2 =1$. 
Let us define $x' := (I + F_{\cJ} C) x$. 
Since $I + F_{\cJ} C$ is assumed to be invertible, one can consider the equivalent system $(I + F_{\cJ})(A+B F_{\cJ} C) (I + F_{\cJ} C)^{-1}  x' = \lambda x'$, $\|x'\|^2=1$. 
By the genericity assumption, the variety $V_{\C}(I)$ has $2 N$ distinct elements with $N =n \cdot 2^{m}$, which proves that the ideal $I$ is zero-dimensional and radical. \\
%
For the last claim, let us first emphasize that the ideal $I$ is not in shape position as for every element $(x, w,\lambda) \in V_{\C}(I)$ one has $(-x,-w,\lambda) \in V_{\C}(I)$. 
Note that the elements of $V_{\C}(I)$ have nonzero $x_1$-coordinates, generically. 
Thus, we consider the change of variables $\hat{x}_1=x_1^2$, $\hat{x}_i = \frac{x_i}{x_1}$, for all $i=2,\dots,n$, $\hat{w}_i = \frac{w_i}{x_1}$, for all $j=1,\dots,m$. 
By homogeneity with respect to $x, w$, one can associate to each $e \in \R[x, w,\lambda]$ a polynomial $\hat{e} \in \R[\hat{x},\hat{w},\lambda]$. 
Namely, one defines 
\[
\hat{e}_i := x_1 e_i = x_1^2 \frac{e_i}{x_1} = (A \hat{x} + B \hat{w} - \lambda \hat{x})_i 
\]
for all $i=1,\dots,n$, 
\[
\hat{e}_{n+i} := x_1^2 \frac{w_i}{x_1} (C \frac{x}{x_1} + D \frac{y}{x_1})-x_1^2 \frac{w_i^2}{x_1^2} = \hat{x}_1 \hat{w}_i (C \hat{x} + D \hat{w})- \hat{x}_1 \hat{w}_i^2,
\]
for all $i=1,\dots,m$, and 
\[
\hat{e}_{n+m+1} = x_1^2 \left(\sum_{i=1}^n \frac{x_i}{x_1^2} + \sum_{i=1}^n \frac{w_i}{x_1^2} \right) - 1 = \hat{x}_1 \left(\sum_{i=1}^n \hat{x}_i^2 + \sum_{i=1}^n \hat{w}_i^2 \right) - 1. 
\]
The ideal $\hat{I}$ generated by $(\hat{e}_i)_{1\leq i \leq n+m+1}$ 
\if{
\begin{align}
\label{eq:idealhat}
\hat{I} := \langle \hat{e}_1, \dots, \hat{e}_{n+m+1}\rangle,
\end{align} 
}\fi
is generically in shape position, thus by Lemma \ref{lm:shape} admits a Gr\"{o}bner basis with the following form: 
\begin{align*}
\{u_0,\hat{x}_1 - u_1, \hat{x}_2 - u_2,  \dots, \hat{x}_n - u_n, \hat{w}_1 - v_1, \dots,   \hat{w}_m - v_m \},
\end{align*}
where $u_i, v_j$ are polynomials in $\R[\lambda]$ and $\deg u_0=\# V_{\C}(\hat{I}) = N$, yielding the desired result.
\end{proofof}
Note that Proposition \ref{prop:zerodim} implies in particular that the feasibility problem \eqref{eqlinear}-\eqref{ineqquadr} has a finite number of solutions, generically. 
A similar structure exploitation has been used in \cite{gradsos} to optimize a given polynomial in $n$ variables over the ambient space $\R^n$, after adding the constraints that all its partial derivatives should cancel. 
\subsection{Complete moment-SOS hierarchies}
\label{sec:complete}
Given $f \in \R[x, w,\lambda]$, we focus on the following polynomial optimization problem (POP): 
\begin{equation}
\label{eq:pop}
f^{\min}: \begin{cases}
\begin{aligned}
\min_{x, w,\lambda} \quad  & f  \\	
\text{s.t.}
\quad & (x, w,\lambda) \text{ satisfies } \eqref{eqlinear}-\eqref{ineqquadr}. 
\end{aligned}
\end{cases}
\end{equation}
%
For all positive integers $n,N$ let $\N_N^n := \{ (\alpha_1,\dots,\alpha_n) \in \N^n : \sum_{i=1}^n \alpha_i \leq N \}$ be the set of multi-indices corresponding to supports of all monomials in $n$ variables with degree at most $N$. 
The dimension of $\N_N^n$ is $\binom{n+N}{n}$. 
We now consider the so-called \textit{moment} hierarchy \cite{lasserre2001global} for the above POP \eqref{eq:pop}, that is the sequence of convex optimization programs, indexed by integers $N \geq 1$:
\begin{equation}
\label{eq:moment}
f_{N}^{\mom}: \begin{cases}
\begin{aligned}
\inf_{y} \quad  & L_y(f)   \\	
\text{s.t.}
\quad & H_N(y) \succeq 0, \\
\quad & H_{N-1}(e y)=0, \quad \forall e \in E, \\
\quad & H_{N-1}(\lambda y) \succeq 0, \quad H_{N-1}(g y) \succeq 0, \quad \forall g \in G, \\
\end{aligned}
\end{cases}
\end{equation}
where 
\begin{itemize}
\item the decision variable vector $y$ is indexed by all elements of $\N_{2 N}^{n+m+1}$;
\item the linear functional $L_y : \R[x, w,\lambda] \to \R$ is defined by $L_y(p)   := \sum_{\alpha} p_{\alpha} y_{\alpha}$, for all $p = \sum_{\alpha} p_\alpha x_1^{\alpha_1} \dots x_n^{\alpha_n} w_1^{\alpha_{n+1}} \dots w_m^{\alpha_{n+m}} \lambda^{\alpha_{n+m+1}}$;
\item the multivariate Hankel -- also called \emph{pseudo-moment} -- matrix $H_N(y)$ is indexed by all multi-indices in $\N_N^{n+m+1}$ and its $(\beta, \gamma)$ entry is $y_{\beta+\gamma}$;
\item the localizing matrix $H_{N-1}(p y)$ associated to a polynomial $p \in \R[x, w,\lambda]$ is indexed by all multi-indices in $\N_{N-1}^{n+m+1}$ and its  $(\beta, \gamma)$ entry is $\sum_{\alpha} p_{\alpha} y_{\alpha+ \beta+\gamma}$. 
\end{itemize}  
The dual of \eqref{eq:moment} is the following sum of squares (SOS) program:
\begin{equation}
\label{eq:sos}
f_{N}^{\sos}: \begin{cases}
\begin{aligned}
\sup \quad  & \rho    \\	
\text{s.t.}
\quad & f - \rho =  \sum_{\ell} u_{0,\ell}^2 + \lambda \sum_{\ell}  u_{\lambda,\ell}^2 +  \sum_{g \in G} g \sum_{\ell} u_{g,\ell}^2 
+ \sum_{e \in E} p_e e, \\
\quad & \rho \in \R, \quad u_{0,\ell}, u_{\lambda,\ell}, u_{g,\ell}, p_e \in \R[x, w,\lambda], \quad \deg u_{0,\ell}, u_{\lambda,\ell} \leq N,   \\
\quad & \deg u_{g,\ell}  \leq N-1, \quad \deg p_e  \leq 2N-2, \quad  \forall  \ell, g \in G, e \in E.    
\end{aligned}
\end{cases}
\end{equation}
We refer the interested reader to \cite{lasserre2001global} for more details on the moment-SOS hierarchy. In particular all the Hankel and localizing matrices are real symmetric matrices and for each $N \geq 1$, each pair  \eqref{eq:moment}-\eqref{eq:sos} can be cast as a primal-dual LMI problem whose optimal value is a lower bound of $f^{\min}$. 
\begin{proposition}
\label{prop:cvg}
Let us assume that the ideal $I$ defined in \eqref{eq:ideal} is zero-dimensional. 
The moment-SOS hierarchy \eqref{eq:moment}-\eqref{eq:sos} converges in finitely many steps to the solution of  POP \eqref{eq:pop}, i.e., $f_N^{\mom} =  f_N^{\sos} = f^{\min}$ for $N$ large enough. 
For $N$ large enough, there is a positive integer $s \leq N$ such that $\rank H_s(y) =  \rank H_{s-1}(y)$ and the set of minimizers of POP \eqref{eq:pop} is equal to $V_{\C} (\kernel H_s(y))$ for any maximal rank optimal solution $y$ of SDP \eqref{eq:moment}. 
\end{proposition}
\begin{proofof}{\rpr{prop:cvg}}
The first claim readily follows from \cite[Theorem~6.15]{laurent2009sums}. 
The second claim is obtained by combining Theorem~6.18 and Theorem~6.20 from \cite{laurent2009sums}.  
\end{proofof}
~\\
An immediate corollary of Proposition \ref{prop:cvg} is that when POP \eqref{eq:pop} has a unique minimizer then for large enough $N$, there is some positive integer $s \leq N$ such that $\rank H_s(y) = 1$ for any maximal rank optimal solution $y$ of SDP \eqref{eq:moment}. 
In this case, the coordinates of the minimizer can be directly extracted from the entries of the vector $y$ indexed by $\alpha$ such that $\sum_{i=1}^{n+m+1}\alpha_i = 1$.  \\
Note that the maximality assumption is not restrictive as most interior-point algorithms used for solving LMI programs return a maximum rank optimal solution; see, e.g.,  \cite[Chapter~4]{de2006aspects} for more details. 
Eventually, if one selects $f = \lambda$ and obtains finite convergence of the hierarchy with $f_N^{\mom} =  f_N^{\sos} = f^{\min} > 0$, then one ensures that the initial problem \eqref{eqlinear}-\eqref{ineqquadr} is feasible with $\lambda > 0$, which proves that the system $\Sigma$ given in \eqref{eq:G} is unstable. 

Let us define $K :=  \R[x, w,\lambda] / I$ and 
\[\B_N:=\{(\lambda^k, \lambda^k x_i, \lambda^k w_j) \mid 0 \leq k \leq N, 1 \leq i \leq n, 1 \leq j \leq m\}, \]
for all $N \in \N$. 
\begin{lemma}
\label{lemma:standard}
Let us assume that the ideal $I$ defined in \eqref{eq:ideal} is as in Proposition \ref{prop:zerodim}.
There exists $N \in \N$ such that the set $\B_N$ includes a basis of $K$. 
\end{lemma}
\begin{proofof}{\rle{lemma:standard}}
By Proposition \ref{prop:zerodim}, the ideal $I$ admits a Gröbner basis with the following form: 
\begin{equation}
\{u_0,x_1^2-u_1,x_2-u_2 x_1,\dots,x_n-u_n x_1,w_1-v_1 x_1,\dots, w_m - v_m x_1\},
\end{equation}
where $u_i,v_j$ are polynomials in $\R[\lambda]$ and $\deg u_0=\# V_{\C}(I)$. 
Let $N = \deg u_0 -1$. 
The set of monomials that cannot be divided by any term from $\LT(I)$ is included in $\B_N$, yielding the desired result. 
\end{proofof}
%
Let $S = \{(x, w,\lambda) \in \R^{n+m+1} : g(x, w,\lambda) \geq 0, e(x, w,\lambda)=0, \forall g \in G, \forall e \in E \}$. 
Given $\ell, N \in \N$, let us denote by $\R[x, w,\lambda]_{=\ell}$ the set of polynomials being homogeneous of degree $\ell$ in $(x, w)$ and $\R[x, w,\lambda]_{N}$ the set of polynomials of total degree less than $d$ in $(x, w, \lambda)$.   
For instance $f = \lambda (\sum_{i=1}^n x_i^2 +  \sum_{j=1}^m w_j^2) \in \R[x, w,\lambda]_{=2}$ and $f \in \R[x, w,\lambda]_{3}$. 
Given two sequences $p$ and $q$ of $n$ polynomials, we note $p \cdot q = \sum_{i=1}^n p_i q_i$. 
\begin{lemma}
\label{lemma:nonneg}
Let us assume that the ideal $I$ defined in \eqref{eq:ideal} is as in Proposition \ref{prop:zerodim} and 
that $f \in \R[x, w,\lambda]_{=2}$ is nonnegative on $S$. 
Then there exist $t, N \in \N$, finite sequences $(u_{0,\ell})_{\ell}$, $(u_{\lambda,\ell})_{\ell}$, $(u_{g,\ell})_{\ell}$, $p=(p_i)_{1\leq i \leq n},  q=(q_j)_{1\leq j \leq m}$ of polynomials such that
\begin{align}
f = & \sum_{\ell} u_{0,\ell}^2 +  \lambda \sum_{\ell} u_{\lambda,\ell}^2 +  \sum_{g \in G} g \sum_{\ell} u_{g,\ell}^2 \\ \nonumber
& + p \cdot (A x + B w - \lambda x) + q \cdot w \odot (C x + D w - w) \\ 
& + r (\|x \|^2 + \|w \|^2 - 1),
\end{align}
with $v_{0,\ell}, v_{\lambda,\ell}, v_{g,\ell}$ supported on $\B_N$, $p_i \in \R[x, w,\lambda]_{2d-1}$, $q_j, r \in  \R[x, w,\lambda]_{2d-2}$, for all $i=1,\dots,n$, $j = 1,\dots,m$, $g \in G$.
\end{lemma}
\begin{proofof}{\rle{lemma:nonneg}}
By Lemma \ref{lemma:standard}, there exists  $N \in \N$ such that the set $\B_N$ contains a  basis of $K$. 
The desired representation then follows by \cite[Theorem 6.15]{laurent2009sums} and its proof (see also \cite[\S~8.2]{laurent2009sums}), where the degree $d$ depends on the degree of $f$ and $N$. 
\end{proofof}
We are now ready to prove our main result. 
Let $\B_N^- := \{(\lambda^k x_i, \lambda^k w_j) \mid 0 \leq k \leq N, 1 \leq i \leq n, 1 \leq j \leq m\}$. 
\begin{theorem}
\label{th:sparse}
Let us assume that the ideal $I$ defined in \eqref{eq:ideal} is as in Proposition \ref{prop:zerodim}. 
Given $f \in \R[x, w,\lambda]_{=2}$, let $f^{\min}$ be its minimum on $S$. 
Let $t, N \in \N$ be as in Lemma \ref{lemma:nonneg}. 
Then there exist finite sequences $(v_{0,\ell})_{\ell}$, $(v_{\lambda,\ell})_{\ell}$, $(v_{g,\ell})_{\ell}$, $p=(p_i)_{1\leq i \leq n},  q=(q_j)_{1\leq j \leq m}$ of polynomials such that
\begin{align}
\label{eq:sparse}
\|(x,w)\|^{2d - 2} (f - f^{\min} \|(x,w)\|^2)  = & \|(x,w)\|^{2d - 2}  \left( \sum_{\ell} v_{0,\ell}^2 + \lambda \sum_{\ell} v_{\lambda,\ell}^2 \right) \nonumber \\  
& + \|(x,w)\|^{2d - 4} \sum_{g \in G} g \sum_{\ell} v_{g,\ell}^2  \\
& + p \cdot (A x + B w - \lambda x) \nonumber\\
& + q \cdot w \odot (C x + D w - w) ,\nonumber
\end{align}
with $v_{0,\ell}, v_{\lambda,\ell}, v_{g,\ell}$ supported on $\B_N^-$, $p_i \in \R[x, w,\lambda]_{=2d-1}$, $q_j \in  \R[x, w,\lambda]_{=2d-2}$, for all $i=1,\dots,n$, $j = 1,\dots,m$, $g \in G$.
\end{theorem}
\begin{proofof}{\rth{th:sparse}}
The polynomial $f - f^{\min} \in  \R[x, w,\lambda]_{=2}$ is nonnegative on $S$ thus by Lemma  \ref{lemma:nonneg} there exist finite sequences $(u_{0,\ell})_{\ell}$, $(u_{\lambda,\ell})_{\ell}$, $(u_{g,\ell})_{\ell}$ of polynomials supported on $\B_N$ such that
\begin{align*}
f - f^{\min} = & \sum_{\ell} u_{0,\ell}^2 + \lambda \sum_{\ell} u_{\lambda,\ell}^2 + \sum_{g \in G} g \sum_{\ell} u_{g,\ell}^2 \\
& + p \cdot (A x + B w - \lambda x) + q \cdot w \odot (C x + D w - w) \\ 
& + r (\|x \|^2 + \|w \|^2 - 1),
\end{align*}
with $p_i \in \R[x, w,\lambda]_{2d-1}$, $q_j, r \in  \R[x, w,\lambda]_{2d-2}$, for all $i=1,\dots,n$, $j = 1,\dots,m$. \\
We shall now rely on a homogenization argument. 
Let us assume that $\|(x, w)\| \neq 0$. 
In the above representation we replace $x$ and $w$ by $x/\|(x, w)\|$ and $w/\|(x, w)\|$, respectively,  and multiply each side by $\|(x, w)\|^{2d}_2$. 
This yields the following representation
\begin{align}
\label{eq:quadproof}
\|(x,w)\|^{2d - 2} (f - f^{\min} \|(x,w)\|^2)  = & \|(x,w)\|^{2d}  \left( \sum_{\ell}  \tilde{u}_{0,\ell}^2 + \lambda \sum_{\ell}  \tilde{u}_{\lambda,\ell}^2 \right) \nonumber \\  
& + \|(x,w)\|^{2d - 2} \sum_{g \in G} g \sum_{\ell}  \tilde{u}_{g,\ell}^2  \\
& + \|(x,w)\|^{2d - 1} \tilde{p} \cdot (A x + B w - \lambda x) \nonumber\\
& + \|(x,w)\|^{2d - 2} \tilde{q} \cdot w \odot (C x + D w - w) ,\nonumber
\end{align}
where $\tilde{h}(x, w,\lambda) :=  h(x/\|(x, w)\|,w/\|(x, w)\|,\lambda)$ for any polynomial $h$. 
Then notice that for each $\ell$, $u_{0,\ell}$ (resp.~$u_{\lambda,\ell}$, $u_{g,\ell}$ for all $g \in G$) is of the form $l_0(\lambda) + l_1(x, w,\lambda)$, with $l_0 \in \R[\lambda]$ and $l_1 \in \R[x, w,\lambda]_{=1}$. \\
Thus, each term $\|(x, w)\|^2_2 \tilde{u}_{0,\ell}^2$ is of the form $\|(x, w)\|^2 \, l_0^2 + l_1^2 + 2 \|(x, w)\| \, l_0 \, l_1$. 
Since the left-hand-side of \eqref{eq:quadproof} is a polynomial, the sum of all terms in the right-hand-side involving odd powers of $\|(x, w)\|$ must cancel. \\
Removing such terms yields a decomposition $ \|(x, w)\|^2 \, l_0^2 + l_1^2 = \sum_{i=1}^n(x_i l_0)^2 +  \sum_{j=1}^m(w_j l_0)^2 + l_1^2 $, where each squared polynomial is supported on $\B_N^-$. 
A similar consideration on the other multipliers yields the desired representation \eqref{eq:sparse} for all $\|(x, w)\| \neq 0$. 
If $\|(x, w)\| = 0$ then the representation is also valid. 
\end{proofof}
The homogenization trick used in the proof of Theorem~\ref{th:sparse} is similar to the one provided in \cite[Theorem~3]{de2005equivalence} to obtain sums of squares decomposition of homogeneous polynomials over the unit sphere. \\
Let us consider the following moment hierarchy doubly indexed by $d,N\in \N$:
\begin{equation}
\label{eq:sparsemoment}
\tilde{f}_{d,N}^{\mom}: \begin{cases}
\begin{aligned}
\inf_{y} \quad  & L_y(\|(x, w)\|^{2d-2}_2 f)   \\	
\text{s.t.}
\quad & \tilde{H}_{d,N}(y) \succeq 0, \\
\quad & \tilde{H}_{d-1,N-1}(\lambda y) \succeq 0, \quad \tilde{H}_{d-2,N}(g y) \succeq 0, \quad \forall g \in G, \\
\quad & L_y(p e_i)=0, \ \forall p \in \R[x, w,\lambda]_{=2d-1}, \deg p \leq 2N+2d-1, \ i=1,\dots,n, \\
\quad & L_y(q e_{n+j})=0, \  \forall q \in \R[x, w,\lambda]_{=2d-2}, \deg q \leq 2N+2d-2, \ j=1,\dots,m, \\
\quad & L_y(\|(x, w)\|^{2d}_2)=1,
\end{aligned}
\end{cases}
\end{equation}
where
\begin{itemize}
\item the multivariate Hankel matrix $\tilde{H}_{d,N}(y)$ is indexed by elements from $\|(x, w)\|^{2d}_2 \B_N^- = \{ \|(x, w)\|^{2d}_2 b : b \in \B_N^-  \}$;
\item the localizing matrix $\tilde{H}_{d-1,N-1}(\lambda y)$ is indexed by elements from $\|(x, w)\|^{2d-2}_2 \B_{N-1}^-$; 
\item the localizing matrix $\tilde{H}_{d-2,N}(g y)$   is indexed by elements from $\|(x, w)\|^{2d-4}_2 \B_{N}^-$. 
\end{itemize}
Considering the dual of \eqref{eq:sparsemoment} yields the following SOS hierarchy:
\begin{equation}
\label{eq:sparsesos}
\tilde{f}_{d,N}^{\sos}: \begin{cases}
\begin{aligned}
\sup_{\rho \in \R} \quad  & \rho &   \\	
\text{s.t.}
\quad & \|(x,w)\|^{2d - 2} (f - \rho \|(x,w)\|^2)  = \|(x,w)\|^{2d - 2}  \left( \sum_{\ell} v_{0,\ell}^2 + \lambda \sum_{\ell} v_{\lambda,\ell}^2 \right)  \\  
\quad & \qquad \qquad \qquad \qquad \qquad \qquad \qquad + \|(x,w)\|^{2d - 4} \sum_{g \in G} g \sum_{\ell} v_{g,\ell}^2  \\
\quad & \qquad \qquad \qquad \qquad \qquad \qquad \qquad  + p \cdot (A x + B w - \lambda x) \\
\quad & \qquad \qquad \qquad \qquad \qquad \qquad \qquad  + q \cdot w \odot (C x + D w - w) ,\\
\quad & v_{0,\ell}, v_{g,\ell} \text{ supported on } \B_N^-, \ v_{\lambda,\ell} \text{ supported on } \B_{N-1}^-, \  \forall  \ell, g \in G,  \\
\quad &  p_i \in \R[x,w,\lambda]_{=2d-1}, \ \deg p_i \leq 2N+2d-1, \ i=1,\dots,n, \\
\quad &  q_j \in \R[x,w,\lambda]_{=2d-2}, \ \deg q_j \leq 2N+2d-2, \ j=1,\dots,m.
\end{aligned}
\end{cases}
\end{equation}
\begin{proposition}
\label{prop:sparse}
Let us assume that the ideal $I$ defined in \eqref{eq:ideal} is as in Proposition \ref{prop:zerodim}. 
Let $f \in \R[x, w,\lambda]_{=2}$ and $f^{\min}$ be its minimum on $S$. Then the following hold: \\
(i) The moment-SOS hierarchy \eqref{eq:sparsemoment}-\eqref{eq:sparsesos} converges in finitely many steps. \\
(ii) For large enough $d,N \in \N$ there exist a positive integer $s \leq N$ such that $\rank \tilde{H}_{d,s}(y) = \rank \tilde{H}_{d,s-1}(y)$. 
%
\end{proposition}
\begin{proofof}{\rpr{prop:sparse}}
(i) As an immediate consequence of Theorem \ref{th:sparse} one obtains finite convergence of the moment-SOS hierarchy \eqref{eq:sparsemoment}-\eqref{eq:sparsesos}. \\
(ii) One can adapt the proof of Theorem~6.20 from \cite{laurent2009sums} to show that for large enough $d,N \in \N$, there is a positive integer $s \leq N$ such that $\rank \tilde{H}_{d,s-1}(y) = \rank \tilde{H}_{d,s}(y)$. We briefly outline the proof strategy: 
\begin{itemize}
\item Let us consider the polynomial $u_0$ of degree $N$ from the set $\Gb$ given in \eqref{eq:standard}, and assume without loss of generality than $u_0 = \lambda^N - \sum_{k=0}^{N-1} u_{0,k} \lambda^k$. 
This polynomial $u_0$ belongs to the ideal generated by the polynomials involved in the equality constraints. 
Thus, a similar homogenization strategy as in the proof of Theorem \ref{th:sparse} yields 
$\|(x,w)\|^{2d} u_0 x_i w_j = p \cdot (A x + B w - \lambda x)  + q \cdot w \odot (C x + D w - w)$, for some finite polynomial sequences $p \subset \R[x,w,\lambda]_{=2d-1}, q \subset \R[x,w,\lambda]_{=2d-2}$ with $d \in \N$ depending on $N$. 
%
\item By using the equality constraints from \eqref{eq:sparsemoment}, we obtain 
\[L\left(\|(x,w)\|^{2d}  \lambda^N x_i w_j \right)= \sum_{k=0}^{N-1} u_{0,k} L \left( \|(x,w)\|^{2d}  \lambda^{k} x_i w_j \right).
\] 
Note that the left-hand side of the above equality is the entry of the matrix  $\tilde{H}_{d,N}(y)$ in the row indexed by $x_i$ and the column indexed by $\lambda^N w_j$. 
The right-hand-side is a weighted sum of entries in the same row and columns indexed by $\lambda^k w_j$, $k=0,\dots,N-1$. 
Similar equalities holds if one replaces  $x_i w_j$ by either $x_i x_j$ or $w_i w_j$. 
Therefore, we can proceed exactly as in the proof of  \cite[Theorem~6.20]{laurent2009sums}: for some large enough $d, N \in \N$, there exists a positive integer $s \leq N$ such that every column of $\tilde{H}_{d,s}(y)$ can be written as a linear combination of columns of $\tilde{H}_{d,s-1}(y)$. 
\end{itemize}
\if{
As in Proposition \ref{prop:cvg} the existence of a rank 1 solution follows by combining Theorem~6.18 and Theorem~6.20 from \cite{laurent2009sums}. 
Let us denote by $\delta_1$ and $\delta_2$ the two Dirac measures supported on $m_1$ and $m_2$, respectively. 
Let $y_1$ and $y_2$ be the moment sequences associated to $\delta_1$ and $\delta_2$, respectively. 
For instance $y_1$ is the infinite sequence $\{ x_1^{\alpha_1} \dots x_n^{\alpha_n} w_1^{\alpha_{n+1}} \dots w_m^{\alpha_{n+m}} \lambda^{\alpha_{n+m+1}} : \alpha \in \N^{n+m+1} \}$. \\
Let us consider the following infinite-dimensional optimization problem over the set $\Pc(S)$ of probability measures supported on $S$:
\begin{align}
\label{eq:meas}
\inf \left\{ \int_S f \text{d} \mu : \mu \in \Pc(S) \right\}.
\end{align}
By a standard argument (see for instance \cite{lasserre2001global}), the optimal value of \eqref{eq:meas} is $f^{\min}$, and the infimum is attained at any convex combination of $\delta_1$ and $\delta_2$. 
We now rely on a sign symmetry argument: given $d \in \N$ and a multi-index $\alpha$ of element from $\B_d^-$, the corresponding moments of $\delta_1$ and $\delta_2$ coincide, as well as the ones of any convex combination of $\delta_1$ and $\delta_2$. 
}\fi
\end{proofof}
We end this section with the following conjecture. 
\begin{conjecture}
\label{conj:rank1}
Let us assume that the ideal $I$ defined in \eqref{eq:ideal} is as in Proposition \ref{prop:zerodim}. 
Let $f \in \R[x, w,\lambda]_{=2}$ and $f^{\min}$ be its minimum on $S$. 
Assume that $f$ has two global minimizers $(\lambda, x, w)$ and $(\lambda, -x, -w)$ on $S$. Then there exist $d,N \in \N$ such that $\rank \tilde{H}_{d,N}(y) = 1$ for any optimal solution $y$ of SDP \eqref{eq:sparsemoment}.
\end{conjecture}
\subsection{Relation with the hierarchy \texorpdfstring{\eqref{eq:dLMIN}}{(26)} from Section \texorpdfstring{\ref{sec:subhigh}}{4.2}}
\label{sec:micnon}
We now relate the SDP hierarchy \eqref{eq:dLMIN} from Section \ref{sec:subhigh} and the above moment hierarchy \eqref{eq:sparsemoment}. \\
In Section \ref{sec:subhigh}, we considered the set $G = \{w_i w_j, w_j w_i - w_j (C x + D w)_i,  [w_i - (C x + D w)_i] [w_j - (C x + D w)_j] \mid 1\leq i,j \leq m \}$ and an objective function that corresponds to a sum of monomial squares $f = \sum_{k=0}^N \lambda^{2k} (\|x\|^2 + \|w\|^2)$. 
\\
By Proposition \eqref{prop:zerodim}, the ideal $I$ defined in \eqref{eq:ideal} is zero-dimensional, generically. 
So the set $S$ is either empty of contains finitely many elements, generically.\\
When $I$ is zero-dimensional and $S$ is non-empty then the minimum $f^{\min}$ of $f$ on $S$ is equal to $\lambda^{\min} := \min \{\lambda : (\lambda, x, w) \in S \}$. 
If $\lambda^{\min} > 0$ then as a direct consequence of Proposition \ref{prop:sparse}, the moment-SOS hierarchy \eqref{eq:sparsemoment}-\eqref{eq:sparsesos} converges in finitely many steps and outputs a proof of instability. \\
The hierarchy \eqref{eq:dLMIN} corresponds to a variant of the above moment hierarchy \eqref{eq:sparsemoment} involving matrix variables of smaller size. 
More specifically, this variant is obtained after selecting $d=1$ in the objective function, $d=0, 1, 2$, in the first, second and third inequality constraint, respectively, as well as $d=1$ in the equality constraints. 
In particular the matrix variable $\tilde{H}_N$ from \eqref{eq:dLMIN} corresponds to $\tilde{H}_{0,N}(y)$. \\
Additionally, one can enforce the polynomial $f - f^{\min} \|x\|^2$ to have a representation as in \eqref{eq:sparse} by imposing that $f^{\min} \|w\|^2$ should be decomposed as $\sum_{g \in G} g \sum_{\ell} v_{g,\ell}^2$ for some $v_{g,\ell} \in \R[\lambda]$. 
Then the last equality constraint of \eqref{eq:sparsemoment} is replaced by $L_y(\|x\|^2) = 1$, which corresponds to the trace equality constraint mentioned at the end of Section \ref{sec:subhigh}. \\
The equality constraints from \eqref{eq:dLMIN} correspond to the first and second set of equality constraints from \eqref{eq:sparsemoment} with $d=1$. 
The set of SDP constraints from  \eqref{eq:dLMIN} correspond to the positive semi-definiteness of principal sub-matrices of $\tilde{H}_{0,N-1}(\lambda y)$. 
The nonnegativity constraints from \eqref{eq:dLMIN} correspond to the nonnegativity of $1\times 1$ principal minors of the localizing matrices $\tilde{H}_{0,N}(g y)$. \\
Even though we cannot guarantee that the more economical LMI hierarchy \eqref{eq:dLMIN} derived in Section \ref{sec:subhigh} has the same convergence properties as the one provided in \eqref{eq:sparsemoment}, it sill performs well in practice, as emphasized by the numerical experiments from Section \ref{sub:num}. 


\section{Conclusion}


In this paper, we addressed the problem of stability analysis for feedback systems involving ReLU nonlinearities.
Our approach focused on a semialgebraic set representation that precisely captures the input-output characteristics of ReLU functions.
Based on this representation, we derived a primal linear matrix inequality (LMI) to serve as a stability certificate, and subsequently developed a dual LMI to certify instability.
Furthermore, we provided a hierarchy of LMI relaxations with convergence guarantees that allows one to certify instability. \\
The first order primal-dual relaxation given in Section \ref{sec:stability} has a clear interpretation, as the dual relaxation provides a sufficient condition for the system $\Sigma$ to be stable. 
For higher order relaxations, we plan to investigate the system theoretic interpretation of the dual program.  \\
We also intend to derive tightness conditions for the first order relaxation, in the spirit of previous studies on robustness certification of single hidden-layer ReLU networks (\cite{Zhang_NIPS2020}). \\
%
%
Eventually, Conjecture \ref{conj:rank1} could be addressed by adapting the proof of Theorem~6.11 (initially presented in \cite[Theorem~12]{schweighofer2005optimization}) or Theorem~6.18 from \cite{laurent2009sums}. 
This happens to be rather non-trivial as the results presented in  \cite{schweighofer2005optimization,laurent2009sums} are obtained for  multivariate Hankel matrices indexed by elements of the canonical basis while in our context these matrices are indexed by monomials that are homogeneous in $(x,w)$. 
The resolution of this conjecture is left for future work. 
\paragraph{Acknowledgements} 
The authors benefited from helpful related discussions with Sami Halaseh, Didier Henrion, Monique Laurent, Bernard Mourrain and Mateusz Skomra. 
Our work was supported by JSPS KAKENHI Grant Number
JP21H01354 and Japan Science and Technology Agency (JST) as
part of Adopting Sustainable Partnerships for Innovative Research
Ecosystem (ASPIRE), Grant Number JPMJAP2402. 
Our work has benefited from the AI Interdisciplinary Institute ANITI. ANITI is funded by the France 2030 program under the Grant agreement n$^{\circ}$ANR-23-IACL-0002. 
Our work was also supported by 
the HORIZON–MSCA-2023-DN-JD of the European Commission under the Grant Agreement No 101120296 (TENORS), as well as by 
the National Research Foundation, Prime Minister's Office, Singapore under its Campus for Research Excellence and  Technological Enterprise (CREATE) programme. 


\begin{thebibliography}{EWM{\etalchar{+}}21}

\bibitem[CLMP20]{chen2020semialgebraic}
T.~Chen, J.~B. Lasserre, V.~Magron, and E.~Pauwels.
\newblock Semialgebraic optimization for lipschitz constants of relu networks.
\newblock {\em Advances in Neural Information Processing Systems},
  33:19189--19200, 2020.

\bibitem[CLO13]{cox2013}
D.~Cox, J.~Little, and D.~OShea.
\newblock {\em Ideals, varieties, and algorithms: an introduction to
  computational algebraic geometry and commutative algebra}.
\newblock Springer Science \& Business Media, 2013.

\bibitem[CTH16]{Carrasco_EJC2016}
J.~Carrasco, M.C. Turner, and W.P. Heath.
\newblock {Z}ames–{F}alb multipliers for absolute stability: From {O}'shea's
  contribution to convex searches.
\newblock {\em European Journal of Control}, 28:1--19, 2016.

\bibitem[D\"10]{Dur_2010}
M.~D\"{u}r.
\newblock Copositive programming - a survey.
\newblock In M.~Diehl, F.~Glineur, E.~Jarlebring, and W.~Michiels, editors,
  {\em Recent Advances in Optimization and Its Applications in Engineering},
  pages 3--20. Springer, 2010.

\bibitem[DK06]{de2006aspects}
E.~De~Klerk.
\newblock {\em Aspects of semidefinite programming: interior point algorithms
  and selected applications}, volume~65.
\newblock Springer Science \& Business Media, 2006.

\bibitem[dKLP05]{de2005equivalence}
E.~de~Klerk, M.~Laurent, and P.~Parrilo.
\newblock On the equivalence of algebraic approaches to the minimization of
  forms on the simplex.
\newblock In {\em Positive Polynomials in Control}, pages 121--132. Springer,
  2005.

\bibitem[Ebi12]{Ebihara_2012}
Y.~Ebihara.
\newblock {\em Systems Control Using LMI (in {J}apanese)}.
\newblock Morikita Syuppan Co., Ltd, Tokyo, 2012.

\bibitem[EOH09]{Ebihara_IEEE2009}
Y.~Ebihara, Y.~Onishi, and T.~Hagiwara.
\newblock Robust performance analysis of uncertain {LTI} systems: Dual {LMI}
  approach and verifications for exactness.
\newblock {\em IEEE Transactions on Automatic Control}, 54(5):938--951, 2009.

\bibitem[EWM{\etalchar{+}}21]{Ebihara_EJC2021}
Y.~Ebihara, H.~Waki, V.~Magron, N.~H.~A. Mai, D.~Peaucelle, and S.~Tarbouriech.
\newblock $l_2$ induced norm analysis of discrete-time {LTI} systems for
  nonnegative input signals and its application to stability analysis of
  recurrent neural networks.
\newblock {\em The 2021 ECC Special Issue of the European Journal of Control},
  62:99--104, 2021.

\bibitem[EXM{\etalchar{+}}24]{Ebihara_arXiv2024a}
Y.~Ebihara, D.~Xin, V.~Magron, D.~Peaucelle, and S.~Tarbouriech.
\newblock Local {L}ipschitz constant computation of {ReLU-FNN}s:\ {U}pper bound
  computation with exactness verification.
\newblock In {\em arxiv:2310.11104 in math.OC, also to appear in Proc. the 22nd
  European Control Conference}, 2024.

\bibitem[FMP22]{Fazlyab_IEEE2022}
M.~Fazlyab, M.~Morari, and G.~J. Pappas.
\newblock Safety verification and robustness analysis of neural networks via
  quadratic constraints and semidefinite programming.
\newblock {\em IEEE Transactions on Automatic Control}, 67(1):1--15, 2022.

\bibitem[FS17]{Fetzler_IFAC2017}
M.~Fetzer and C.~W. Scherer.
\newblock Absolute stability analysis of discrete time feedback
  interconnections.
\newblock {\em IFAC PapersOnline}, 50(1):8447--8453, 2017.

\bibitem[GM89]{gianni89}
P.~Gianni and T.~Mora.
\newblock {Algebraic solution of systems of polynomial equations using Gröbner
  bases}.
\newblock In {\em Applied Algebra, Algebraic Algorithms and Error Correcting
  Codes, Proceedings of AAECC-5, volume 356 of LNCS}, pages 247--257. Springer,
  1989.

\bibitem[GR22]{Gronqvist_MTNS2022}
J.~Gr\"{o}nqvist and A.~Rantzer.
\newblock Dissipativity in analysis of neural networks.
\newblock In {\em Proc. the 25th International Symposium on Mathematical Theory
  of Networks and Systems}, pages 1221--1224, 2022.

\bibitem[GVdS19]{Groff_CDC2019}
L.~B. Groff, G.~Valmorbida, and J.~M.~Gomes da~Silva.
\newblock Stability analysis of piecewise affine discrete-time systems.
\newblock In {\em Proc. Conference on Decision and Control}, pages 8172--8177,
  2019.

\bibitem[HL05]{Henrion_2005b}
D.~Henrion and J.~B. Lasserre.
\newblock Detecting global optimality and extracting solutions in gloptipoly.
\newblock In D.~Henrion and A.~Garulli, editors, {\em Lecture Notes in Control
  and Information Sciences}, volume 312. Springer Verlag, Berlin, 2005.

\bibitem[Kha02]{Khalil_2002}
H.~Khalil.
\newblock {\em Nonlinear Systems}.
\newblock Prentice Hall, 2002.

\bibitem[Kor22]{korda2022stability}
M.~Korda.
\newblock Stability and performance verification of dynamical systems
  controlled by neural networks: algorithms and complexity.
\newblock {\em IEEE Control Systems Letters}, 6:3265--3270, 2022.

\bibitem[Las01]{lasserre2001global}
J.~B. Lasserre.
\newblock Global optimization with polynomials and the problem of moments.
\newblock {\em SIAM Journal on optimization}, 11(3):796--817, 2001.

\bibitem[Lau09]{laurent2009sums}
M.~Laurent.
\newblock Sums of squares, moment matrices and optimization over polynomials.
\newblock In {\em Emerging applications of algebraic geometry}, pages 157--270.
  Springer, 2009.

\bibitem[MDV23]{gradsos}
V.~Magron, M.~Safey~El Din, and T.-H. Vu.
\newblock {Sum of Squares Decompositions of Polynomials over their Gradient
  Ideals with Rational Coefficients}.
\newblock {\em SIAM Journal on Optimization}, 33(1):63--88, 2023.

\bibitem[MR97]{Megretski_IEEE1997}
A.~Megretski and A.~Rantzer.
\newblock System analysis via integral quadratic constraints.
\newblock {\em IEEE Transactions on Automatic Control}, 42(6):819--830, 1997.

\bibitem[MS09]{Masubuchi_SCL2009}
I.~Masubuchi and C.~W. Scherer.
\newblock A recursive algorithm of exactness verification of relaxations for
  robust {SDP}s.
\newblock {\em Systems and Control Letters}, 58(8):592--601, 2009.

\bibitem[MW23]{Magron_2023}
V.~Magron and J.~Wang.
\newblock {\em Sparse Polynomial Optimization: Theory and Practice}.
\newblock World Scientific, 2023.

\bibitem[NSY{\etalchar{+}}24]{relufeedback}
S.~Nishinaka, R.~Saeki, T.~Yuno, Y.~Ebihara, V.~Magron, D.~Peaucelle,
  S.~Zoboli, and S.~Tarbouriech.
\newblock Stability analysis of feedback systems with relu nonlinearities via
  semialgebraic set representation.
\newblock {\em {Proceedings of the 4th IFAC Conference of Modelling,
  Identification and Control of nonlinear systems (MICNON)}}, 2024.

\bibitem[O'S67]{O'Shea_IEEE1967}
R.~O'Shea.
\newblock An improved frequency time domain stability criterion for autonomous
  continuous systems.
\newblock {\em IEEE Transactions on Automatic Control}, 12(6):725--731, 1967.

\bibitem[RSL18]{Raghunathan_NIPS2018}
A.~Raghunathan, J.~Steinhardt, and P.~Liang.
\newblock Semidefinite relaxations for certifying robustness to adversarial
  examples.
\newblock {\em Advances in Neural Information Processing Systems}, pages
  10900--10910, 2018.

\bibitem[RWM21]{Revay_LCSS2021}
M.~Revay, R.~Wang, and I.~R. Manchester.
\newblock A convex parameterization of robust recurrent neural networks.
\newblock {\em IEEE Control Systems Letters}, 5(4):1363--1368, 2021.

\bibitem[Sch05a]{Scherer_SIAM2005}
C.~W. Scherer.
\newblock Relaxations for robust linear matrix inequality problems with
  verifications for exactness.
\newblock {\em SIAM Journal on Matrix Analysis and Applications},
  27(2):365--395, 2005.

\bibitem[Sch05b]{schweighofer2005optimization}
Markus Schweighofer.
\newblock Optimization of polynomials on compact semialgebraic sets.
\newblock {\em SIAM Journal on Optimization}, 15(3):805--825, 2005.

\bibitem[Sch06]{Scherer_EJC2006}
C.~W. Scherer.
\newblock {LMI} relaxations in robust control.
\newblock {\em European Journal of Control}, 12(1):3--29, 2006.

\bibitem[Sch22]{Scherer_IEEEMag2022}
C.~W. Scherer.
\newblock Dissipativity and integral quadratic constraints: Tailored
  computational robustness tests for complex interconnections.
\newblock {\em IEEE Control Systems Magazine}, 42(3):115--139, 2022.

\bibitem[TdSJQ11]{Tarbouriech_2011}
S.~Tarbouriech, J.~M.~Gomes da~Silva~Jr, and I.~Queinnec.
\newblock {\em Stability and Stabilization of Linear Systems with Saturating
  Actuators}.
\newblock Springer, 2011.

\bibitem[YSA22]{Yin_IEEE2022}
H.~Yin, P.~Seiler, and M.~Arcak.
\newblock Stability analysis using quadratic constraints for systems with
  neural network controllers.
\newblock {\em IEEE Transactions on Automatic Control}, 67(4):1980--1987, 2022.

\bibitem[ZF68]{Zames_SIAM1968}
G.~Zames and P.~Falb.
\newblock Stability conditions for systems with monotone and slope-restricted
  nonlinearities.
\newblock {\em SIAM Journal on Control}, 6(1):89--108, 1968.

\bibitem[Zha20]{Zhang_NIPS2020}
R.~Y. Zhang.
\newblock On the tightness of semidefinite relaxations for certifying
  robustness to adversarial examples.
\newblock {\em Advances in Neural Information Processing Systems},
  33:3808--3820, 2020.

\end{thebibliography}

\newcommand{\etalchar}[1]{$^{#1}$}

\end{document}